\documentclass[11pt,oneside,reqno]{amsart}
\usepackage{latexsym}
\usepackage{amssymb}
\usepackage{amsthm}
\usepackage{eucal}
\usepackage{amsmath}
\usepackage{amssymb}

\newtheorem{theoremm}{Theorem}
\newtheorem{theorem}{Theorem}[section]
\newtheorem{lemma}[theorem]{Lemma}
\newtheorem{corollary}[theorem]{Corollary}
\newtheorem{proposition}[theorem]{Proposition}

\theoremstyle{definition}

\newtheorem{question}{Question}
\newtheorem{problem}{Problem}
\theoremstyle{remark}

\renewcommand{\bar}{\overline}
\def\supp{\mathop{\rm supp}\nolimits}

\newcommand{\C}{{\mathbb C}}
\newcommand{\N}{{\mathbb N}}

\newcommand{\T}{\mathbb T}
\newcommand{\Z}{{\mathbb Z}}

\newcommand{\cont}{\mathfrak c}

\def\hull#1{\langle#1\rangle}

\author{Jorge Galindo$^{(1)}$ \and Luis Recoder-N\'u\~nez$^{(2)}$ \and
Mikhail Tkachenko$^{(3)}$}

\thanks{This work was started while the second and third
listed authors were visiting the Department of Mathematics
of Universitat Jaume I, Spain, in 2008. They thank their hosts
for hospitality and support.}
\thanks{$^{(1)}$\, Research of this author  supported by the Spanish
Ministry of Science (including FEDER funds), grant MTM2008-04599/MTM
and  Fundaci\'o Caixa Castell\'o-Bancaixa, grant P1.1B2008-26.}
\thanks{$^{(2)}$\, Research of this author supported by a grant
for foreign visitors under the Promotion of Research Program
of University Jaume I}
\thanks{$^{(3)}$\,This author was supported by
CONACyT of Mexico, grant 000000000074468.}
\date{February 15, 2010}

\address{\noindent Jorge Galindo, Instituto de Matem\'aticas y
Aplicaciones (IMAC)\\ Departamento de  Matem\'aticas\\ Universidad
Jaume I, E-12071, Cas\-tell\'on, Spain. \hfill\break \noindent
E-mail: {\tt jgalindo@mat.uji.es}}
\address{\enlargethispage{3cm}\noindent Luis Recoder-N\'u\~{n}ez\\
Department of Mathematical Sciences\\ Central Connecticut State
University\\ 1615 Stanley Street, New Britain, CT 06050,
USA.\hfill\break \noindent \noindent E-mail: {\tt
recoderl@ccsu.edu}}
\address{\enlargethispage{3cm}\noindent  Mikhail Tkachenko
\\Departamento de Matem\'aticas,
Universidad Aut\'onoma Met\-ropolitana\\ Av. San Rafael Atlixco 186,
Col. Vicentina, Iztapalapa, C.P. 09340, M\'exico, D.F.,
Mexico.\hfill\break \noindent E-mail: {\tt mich@xanum.uam.mx}}

\title[Nondiscrete $P$-groups can be reflexive]{Nondiscrete
$P$-groups can be reflexive}

\begin{document}

\begin{abstract}
We present a series of examples of nondiscrete reflexive $P$-groups
(i.e., groups in which all $G_\delta$-sets are open) as well as
noncompact reflexive $\omega$-bounded groups (in which the closure
of every countable set is compact). Our main result implies that every
product of feathered (equivalently, almost metrizable) Abelian groups
equipped with the $P$-modified topology is a reflexive group. In particular,
every compact Abelian group with the $P$-modified topology is reflexive.
This answers a question posed by S.~Hern\'andez and P.~Nickolas
and solves a problem raised by Ardanza-Trevijano, Chasco,
Dom\'{\i}nguez, and Tkachenko.
\end{abstract}

\maketitle

\textit{AMS Subj. Class. (2000):}
\small{54H11; 22D35; 54G10; 54D30; 54B10}


%%%%%%%%%%%%%%%%%%%%%%%%%%%%%%%%%%
\section{Introduction}%%%%%%%%%%%%%%%%%%%%%%%
Extending Pontryagin's duality to diverse classes of topological groups
beyond locally compact ones has been the object of attention through
the last 60 years. It has become patent in recent times that the duality
properties of precompact groups and of projective limits of discrete
groups, two otherwise well studied classes of  topological groups,
are poorly understood. We refer the reader to \cite{chasmart08} for the
case of precompact groups and to \cite{hofmmorr07} and \cite{hofmmorr09}
where the need of an improved knowledge of the duality properties  emerged
in the study of the structure of projective limits of Lie groups.

The duality theory of $P$-groups involves immediately both classes of
topological groups. A $P$-group is a topological group in which the intersection
of any countable family of open sets is open. Some basic information
about $P$-groups can be found in \cite[Section~4.4]{AT}.

It is easy to see that $P$-groups are projective limits of discrete
groups that do not contain infinite compact subsets. Because of this
latter property the character group of a $P$-group is always
precompact. Even more, it has the far stronger property of being
\textit{$\omega$-bounded\/} (the closure of any countable subset is compact),
see Lemma~\ref{omeb} below. Therefore $P$-groups occupy the region
of pro-Lie groups that is farthest to locally compact groups whereas
$\omega$-bounded groups are among the most compact-like classes of
topological groups.

Despite a considerable interest to the subject, the only result in
the literature concerning the duality of $P$-groups is Leptin's
example \cite{lept}, later reproduced by Noble \cite{noble70} and
Banaszczyk \cite[Example~7.11]{bana91}. This example is the inverse
limit of an uncountable family of discrete groups which turns out to
be a nondiscrete $P$-group whose second dual is discrete. Hence
Leptin's group is not reflexive. This motivates the following
question posed by S.~Hern\'andez and P.~Nickolas at the
International Workshop on Topological Groups and Dynamic Systems
in Madrid, 2008:

\begin{question}\label{quest:P}
Must every  reflexive Abelian $P$-group be discrete?
\end{question}

This problem is naturally linked with the question on whether a
precompact, noncompact group can be reflexive (recall that the dual
group of a $P$-group is $\omega$-bounded and hence precompact). The
first examples of precompact, noncompact reflexive groups have  been
recently obtained in \cite{ardaetal} and \cite{galimaca09}, but the
examples presented there are precompact groups with no infinite
compact subsets; so they are far from being $\omega$-bounded.
Question~\ref{quest:P} has therefore the natural accompanying
question.

\begin{question}\label{quest:w}
How close to being compact can a reflexive precompact  Abelian
group be? In particular, can a noncompact $\omega$-bounded group
be reflexive?
\end{question}

The main objective of this paper is to answer Question~\ref{quest:P}
(in the negative) and Question~\ref{quest:w} (in the positive). We
do this in Theorem~\ref{cor:main} below by proving that every product
of discrete Abelian groups with the $P$-modified topology is reflexive.
This result is then extended to products of feathered (equivalently, almost
metrizable) Abelian groups. In an attempt to trace the borders between
reflexive and nonreflexive $P$-groups we also give a number of new
examples of reflexive and nonreflexive $P$-groups.

We also establish in Propositions~\ref{Pro:Une} and~\ref{Pro:Quo} that
the class of reflexive $P$-groups has unexpectedly good permanence
properties---it contains quotient groups and if a reflexive $P$-group $G$
is a dense subgroup of a topological group $H$, then $H$ is reflexive
$P$-group as well.

%%%%%%%%%%%%%%%%%%%%%%%%%%%%%%%%%%%
\section{Notation}\label{Sec:Not}%%%%%%%%%%%%%%%%%%%
All groups considered here are assumed to be Abelian if otherwise is not
specified explicitely. The complex plane with its usual multiplication and
topology is denoted by $\C$. A \emph{character\/} on a group $G$ is a
homomorphism of $G$ to the circle group $\T=\{z\in\C: |z|=1\}$. The set
$\{e^{i\varphi}: -\pi/2<\varphi<\pi/2\}$ is denoted by $\T_+$. Pontryagin's
duality theory is based on relating a topological group $G$ with the group
$G^\wedge$ of continuous characters of $G$. The group $G^\wedge$
will be equipped with the topology of uniform convergence on the compact
subsets of $G$. This topology has as a neighbourhood basis at the identity
the sets
\[
K^\vartriangleright =\left\{ \chi \in G^\wedge \colon
\chi(x)\in\T_+ \mbox{ for all } x\in K\right\},
\]
where $K$ runs over the family of compact subsets of $G$. If a topological
group $G$ is well-represented by its continuous characters, one can
recover $G$ from $G^\wedge$ by forming the bidual group
$G^{\wedge\wedge}=\left(G^\wedge\right)^\wedge$ and then considering
the canonical \textit{evaluation homomorphism\/} $\alpha_{_{G}}\colon G
\to G^{\wedge\wedge}$ defined by $\alpha_{_{G}}(x)(\chi)=\chi(x)$, for all
$x\in G$ and $\chi\in G^\wedge$. We say accordingly that $G$ is
\emph{reflexive} if the homomorphism $\alpha_{_{G}}\colon G\to
G^{\wedge\wedge}$ is a topological isomorphism.

We will use the concept of \textit{nuclear group\/} as it appears in
\cite{Aub,bana91}. It is worth mentioning that every closed subgroup
$H$ of a nuclear group $G$ is \emph{dually embedded\/} in $G$, i.e.,
for every $x\in G\setminus H$ there exists a continuous character
$\chi\colon G\to\T$ such that $\chi(H)=\{1\}$ and $\chi(x)\neq 1$
(see \cite[Corollary~8.6]{bana91}). In particular, continuous
characters of a nuclear group $G$ separate points of $G$ which in
its turn implies that the evaluation homomorphism $\alpha_{_{G}}
\colon G\to G^{\wedge\wedge}$ is injective.

Given a topological group $H$, we denote by $PH$ the
\textit{$P$-modification\/} of $H$ which is the same underlying
group $H$ endowed with the finer topological group topology whose
base is formed by the family of $G_\delta$-sets in the original
group $H$. The subgroup of $H$ generated by a set $A\subseteq H$
is $\hull{A}$. Sometimes we use $\hull{A,B}$ for $A,B\subseteq H$ to
denote the group $\hull{A\cup B}$.

A group $G$ is \emph{boolean\/} if every $x\in G$ satisfies
$x+x=0_G$. The group $\Z_2=\{0,1\}$ and all powers of $\Z_2$ are
boolean groups.

Let $\{D_i\colon i\in I\}$ be a family of topological groups and
$D=\prod_{i\in I}D_i$ the product group with the usual Tychonoff
product topology. Elements of $D$ will be regarded as functions
$x\colon I\to \bigcup_{i\in I}D_i$ such that $x(i)\in D_i$ for all $i\in I$.
If $x\in D$, we define the \textit{support\/} of $x$ as
\[
\supp(x)=\{i\in I\colon x(i)\neq 0_i\},
\]
where $0_i$ is the neutral element of $D_i$. With these notations,
the subgroup
\[
\Sigma{D}=\left\{x\in D\colon |\supp(x)|\leq \omega\right\}
\]
of $D$ is called the \emph{$\Sigma$-product\/} of the family
$\{D_i: i\in I\}$. Similarly,
\[
\sigma{D}=\left\{ x\in D\colon |\supp(x)|< \omega\right\}
\]
is a subgroup of $D$ which is called the \emph{$\sigma$-product\/}
of the family $\{D_i: i\in I\}$.

Let $[I]^{\leq \omega}$ denote the family of all countable subsets of
the index set $I$. For every $i\in I$, we denote by $\mathcal{S}_i$
the family of all subgroups of type $G_\delta$ in $D_i$. It is clear
that $\mathcal{S}_i$ is a base of $PD_i$ at $0_i$. The collection
of sets
\[
U(J,f)=\left\{ x\in D\colon x(i)\in f(i) \mbox{ for all } i\in J\right\},
\]
where $J\in [I]^{\leq \omega}$ and $f$ a function with the domain
$J$ such that $f(i)\in\mathcal{S}_i$ for each $i\in J$, constitutes
a base at the identity of $PD$, the $P$-modification of the product
group $D$. Clearly, $\Sigma{D}$ is a dense subgroup of $PD$,
while $\sigma{D}$ is a dense subgroup of $D$.

In what follows we will also use the sets
\[
U(J)=\left\{ x\in D\colon x(i)=0_i \mbox{ for all } i\in J\right\},
\]
with $J\subset I$. Notice that if the groups $D_i$'s are discrete,
then the family $\{U(J): J\in [I]^{\leq \omega}\}$ forms a local base
at the identity of $PD$.

We will mainly be working with topological groups $G$ such that
\[
\Sigma{D} \subseteq G \subseteq PD.
\]
The subgroup $\Sigma{D}$ of $PD$ will always carry the topology
inherited from $PD$, i.e., $\Sigma{D}$ is a $P$-group.

Let $G$ be a subgroup of $PD$ and $J\subset I$. We will say
that a character $\chi\colon G\to\T$ \emph{depends on $J$} if there is
$x\in G$ with $\supp(x)\subset J$ such that $\chi(x)\neq 1$ (notice that
$x\in G\cap U(I\setminus J)$).

It is easy to see (see Lemma~\ref{Le:Depen} below) that a character
$\chi\colon G\to \T$ is continuous if and only if there are $J\in [I]^{\leq
\omega}$ and a function $f$ such that $G\cap U(J,f)\subset\ker\chi$.
Therefore, for every continuous character $\chi\colon G \to \T$,
there is $J\in [I]^{\leq\omega}$ such that $\chi$ does not depend on
$I\setminus J$. We say in this case that \emph{$\chi$ depends
(at most) on countably many coordinates}. Similarly, we say that a set
\emph{$K\subset G^\wedge$ depends (at most) on countably many
coordinates\/} if there is $J\in [I]^{\leq \omega}$ such that every $\chi\in K$
does not depend on $I\setminus J$.

%%%%%%%%%%%%%%%%%%%%%%%%%%%%%%%%%%
\section{The results}\label{Sec:MR}%%%%%%%%%%%%%%%%
As we mentioned above, Leptin~\cite{lept} gave an example of a
nonreflexive $P$-group which was the subgroup of $P\Z_2^{\omega_1}$
consisting of elements with finite support (i.e. the $\sigma$-product of
$\omega_1$ copies of the discrete two-element group $\Z_2$). Here
we extend Leptin's argument to deduce the following result (see
Proposition~\ref{sigmanot}):

\begin{theoremm}\label{sigma}
The $\Sigma$-product $\Sigma{D}\subset PD$ is not reflexive,
where $D=\prod_{i\in I}D_i$ is the product of an uncountable family
of nontrivial discrete Abelian groups.
\end{theoremm}

We also extend Theorem~\ref{sigma} to certain subgroups between $\Sigma=
\Sigma\Z_2^\tau$ and $P\Z_2^\tau$ as follows (see Theorem~\ref{GLnot}):

\begin{theoremm}\label{sigmal}
If $L$ is a countable subgroup of $P\Z_2^\tau$, for an uncountable
cardinal $\tau$, then the subgroup $G_L=\Sigma+L$ of $P\Z_2^\tau$
is not reflexive.
\end{theoremm}

It may be worth to observe that while the proof of Theorem~\ref{sigma}
uses an argument close to Leptin's, this argument does not work for
$G_L$ and a different one is needed for Theorem~\ref{sigmal}.

Somewhat surprisingly we also find a wealth of non-discrete reflexive
$P$-groups (and hence reflexive noncompact $\omega$-bounded groups).
In particular, we prove the following fact in Theorem~\ref{cor:main} (and
Corollary~\ref{Cor:Exis}):

\begin{theoremm}\label{prod}
Let $D=\prod_{i\in I}D_i$ be a product of nontrivial discrete Abelian groups,
where $|I|>\omega$. Then the nondiscrete $P$-group $\Pi=PD$ and the
noncompact $\omega$-bounded group $\Pi^\wedge$ are reflexive.
\end{theoremm}

The non-reflexive groups in Theorems~\ref{sigma} and~\ref{sigmal}
are evidently non-complete, while the reflexive groups presented in
Theorem~\ref{prod} are complete (apply \cite[Theorem~8]{HT}).
One can conjecture, therefore, that reflexive $P$-groups are complete.
We show in Theorem~\ref{T:SSigma} below that this is not the case.

Theorem~\ref{prod} can be given a considerably more general form.
We recall that a topological group $G$ is called \textit{feathered\/} or
\textit{almost metrizable\/} provided there exists a nonempty compact
set $K$ in $G$ with a countable neighbourhood base. The next result
follows from Theorem~\ref{Th:Gen}:

\begin{theoremm}\label{prodF}
Let $D=\prod_{i\in I}D_i$ be a product of feathered Abelian groups.
Then the group $PD$ is reflexive.
\end{theoremm}

As a step towards the proof of Theorem~\ref{prodF} we show that
every compact Abelian group with the $P$-modified topology is reflexive.

The following reflexion principle turns out to be quite useful
when trying to extend the class of reflexive $P$-groups (see
Theorem~\ref{Reduc}):

\begin{theoremm}\label{T:Reduc}
Suppose that $D=\prod_{i\in I}D_i$ is a product of topological groups
and $\Sigma{D}\subset G\subset PD$. Then the group $G$ is reflexive
iff the subgroup $\pi_J(G)$ of $PD_J$ is reflexive, for every set
$J\subset I$ satisfying $|J|\leq\aleph_1$ (here $\pi_J\colon D\to
\prod_{i\in J} D_i$ is the projection).
\end{theoremm}

Theorem~\ref{T:Reduc} enables us to present examples of dense reflexive
subgroups $G$ of the groups $PD$ (see Corollary~\ref{Cor:Sigma}).
Clearly, $G$ is not complete provided that $G\neq PD$:

\begin{theoremm}\label{T:SSigma}
Let $D=\prod_{i\in I} D_i$ be a product of feathered Abelian groups.
Then the subgroup
\[
\Sigma_{\aleph_1}{D}=\{x\in PD: |\supp(x)|\leq\aleph_1\}
\]
of $PD$ is reflexive and every subgroup $G$ of $PD$
containing $\Sigma_{\aleph_1}{D}$ is also reflexive.
\end{theoremm}

Once we have established that the class of reflexive $P$-groups
is quite wide, it is natural to clarify the permanence properties of
this class. In the next result we present two of them (see
Propositions~\ref{Pro:Une} and~\ref{Pro:Quo}):

 \begin{theoremm}\label{Cor:SSigma}
Let $G$ be a reflexive $P$-group.
\begin{enumerate}
\item[(a)] If $G$ is a dense subgroup of a topological group $H$,
               then $H$ is also a reflexive $P$-group.
\item[(b)] If $\pi\colon G\to K$ is a continuous open homomorphism
               of $G$ onto $K$, then $K$ is a reflexive $P$-group.
\end{enumerate}
\end{theoremm}

Let $\tau>\omega$ be a cardinal. To extend Theorem~\ref{prod} to
subgroups slightly smaller than $\Pi=P\Z_2^\tau$ we consider an
arbitrary ultrafilter $\xi$ on $\tau$ containing all subsets of $\tau$
with countable complement. Note that $A\in\xi$ implies
$|A|>\omega$. We then define
\[
G_\xi=\{x\in \Pi\colon \supp(x)\notin \xi\}.
\]
It is straightforward to check that $G_\xi$ is a subgroup of $\Pi$.
Also, if $x\in\Sigma$, then $|\supp(x)|\leq \omega$. So,
$\supp(x)\notin \xi$ and therefore $\Sigma\Pi\subset G_\xi$. It is
clear that $G_\xi\neq \Pi$ because the constant function
$\mathbf{1}$ is not in $G_\xi$. Actually this is the ``most
important" element absent in $G_\xi$, as $\Pi=\langle G_\xi,
\mathbf{1}\rangle=G_\xi \oplus \langle \mathbf{1}\rangle$. We
prove in Theorem~\ref{FinSec} that this smaller subgroup
$G_\xi$ of $\Pi$ is also reflexive:

\begin{theoremm}
The nondiscrete $P$-group $G_\xi$ is reflexive, for every ultrafilter
$\xi$ on $\tau$ containing the complements to countable sets.
\end{theoremm}

%%%%%%%%%%%%%%%%%%%%%%%%%%%%%%%%%%%%
\section{Factorization of characters on subgroups of a product group}%%
\label{Sec:Fac}%%%%%%%%%%%%%%%%%%%%%%%%%%%%%
Here we collect several results of technical nature which will be used later.
Throughout this section $D=\prod_{i\in I} D_i$ stands for the Tychonoff
product of a family $\{D_i: i\in I\}$ of topological groups, not necessarily
Abelian. For $J\subset I$, we denote by $D_J$ the corresponding
subproduct $\prod_{i\in J}D_i$ and by $\pi_J \colon D\to D_J$ the projection.

\begin{lemma}\label{Le:Depen}
Let $G$ be a subgroup of $PD$ and $\chi\colon G\to\T$ be a
character. The following assertions are then equivalent:
\begin{enumerate}
  \item $\chi$ is continuous, i.e., $\chi \in G^\wedge$.
  \item There are $J\in [I]^{\leq\omega}$ and a function $f$ such that
  $G\cap U(J,f)\subset \ker \chi$ (in particular, $\chi$ does not depend
  on $I\setminus J$).
  \item There are a countable set $J\subset I$ and a continuous character
  $\chi_J$ on the subgroup $\pi_J(G)$ of $PD_J$ such that $\chi=
  \chi_J\circ\pi_J\hskip-3.5pt\restriction_G$.
\end{enumerate}
\end{lemma}

\begin{proof}
Suppose $\chi$ is continuous. Since $G$ is a $P$-group, the kernel
of $\chi$ is an open subgroup of $G$. It follows from the definition
of the topology of $PD$ that there exist a countable set $J\subset I$
and a function $f$ with the domain $J$ such that $f(i)$ is an open
subgroup of $PD_i$ for each $i\in J$ and the basic open set
$U(J,f)= \left\{ x\in D\colon x(i)\in f(i) \mbox{ for all } i\in J\right\}$ in
$PD$ satisfies $G\cap U(J,f)\subset\ker\chi$. Clearly $\chi$ does not
depend on $I\setminus J$. Hence (1) implies (2).

Suppose now that $G\cap U(J,f)\subset \ker\chi$ for a countable
set $J\subset I$ and a corresponding function $f$. We define a character
$\chi_J$ on $\pi_J(G)$ by $\chi_J(\pi_J(x))=\chi(x)$ for any element $x\in G$.
This definition is correct since the equality $\pi_J(x)=\pi_J(y)$ implies
that $x^{-1}y\in G\cap U(J,f)$ and $\chi(x)=\chi(y)$. By the definition of
$\chi_J$, we see that $\chi=\chi_J \circ\pi_J\hskip-3.5pt\restriction_G$.
It also follows from the definition of $\chi_J$ that its kernel contains the
set $\pi_J(G)\cap U_f$, where $U_f=\prod_{i\in J} f(i)$ and $f(i)$ is an open
subgroup of $PD_i$ for each $i\in J$. Since $U_f$ is an open subgroup
of $PD_J$, we conclude that $\chi_J$ is a continuous character on the group
$\pi_J(G)$. So (2) implies (3).

Finally suppose that  there is a countable set $J\subset I$ and a continuous
character $\chi_J$ on $\pi_J(G)$ such that $\chi=\chi_J\circ\pi_J\hskip-3.5pt
\restriction_G$. Since the projection $\pi_J\colon PD\to PD_J$ is continuous,
we see that so is the character $\chi$. Hence (3) implies (1).
\end{proof}

The following result is close in the spirit to \cite[Theorem~4.6]{comgore09},
where the product space carries the usual Tychonoff product topology.

\begin{lemma}\label{Le:Dep3}
Let $G$ be a dense subgroup of $PD$ and $\chi\in G^\wedge$.
Then $\chi$ admits a continuous extension to a character
$\bar{\chi}$ on $PD$ and, for every set $J\subset I$, $\bar{\chi}$
does not depend on $I\setminus J$ if and only if there exists a
continuous character $\psi$ on $\pi_J(G)$ such that
$\chi=\psi\circ\pi_J\hskip-4pt\restriction_G$.
\end{lemma}

\begin{proof}
Since $G$ is dense in $PD$ and the circle group $\T$ is compact
(hence complete), $\chi$ extends to a continuous character
$\bar\chi$ on the group $PD$. If $\bar\chi$ does not depend on
$I\setminus J$, where $J\subset I$, there exists a character
$\bar\psi$ on $PD_J$ such that $\bar\chi=\bar\psi\circ\pi_J$,
where $D_J=\prod_{i\in I}D_i$ and $\pi_J\colon D\to D_J$ is the
projection. Since $\pi_J\colon PD\to PD_J$ is open, the character
$\bar\psi$ on $PD_J$ is continuous. Then $\chi=\psi\circ\pi_J
\hskip-4pt\restriction_G$, where $\psi$ is the restriction of $\bar\psi$
to the subgroup $\pi_J(G)$ of $PD_J$. Clearly, the character $\psi$
is continuous.

Conversely, suppose that there exists a continuous character
$\psi$ on the subgroup $\pi_J(G)$ of $PD_J$ such that $\chi=
\psi\circ\pi_J\hskip-3.5pt\restriction_G$. Since $\pi_J(G)$ is dense
in $PD_J$, $\psi$ extends to a continuous character $\bar\psi$ on
the group $PD_J$. Clearly, the characters $\bar\chi$ and
$\bar\psi\circ\pi_J$ coincide on the dense subgroup $G$ of $PD$.
Since the group $\T$ is Hausdorff, we see that $\bar\chi=
\bar\psi\circ\pi_J$. It follows that $\bar\chi$ does not depend on
$I\setminus J$.
\end{proof}

In general, the existence of a continuous character $\psi$ on
$\pi_J(G)$ satisfying $\chi=\psi\circ\pi_J\hskip-4pt\restriction_G$
in Lemma~\ref{Le:Dep3} cannot be weakened to a simpler
condition that $\chi$ does not depend on $I\setminus J$ . However,
this weakening is possible for special subgroups of $PD$ as we will
see in Corollary~\ref{Cor:N} below. First we need a lemma.

\begin{lemma}\label{Le:Open}
If $\,\Sigma{D}\subset G\subset PD$, then the restriction of the
projection $\pi_J\colon PD\to PD_J$ to $G$ is an open homomorphism
of $G$ onto $\pi_J(G)$, for every nonempty set $J\subset I$.
\end{lemma}

\begin{proof}
Let $J$ be a nonempty subset of $I$. Since $G$ is a subgroup of
$PD$, the restriction of $\pi_J$ to $G$ is a continuous homomorphism.
Hence it suffices to verify that the image $\pi_J(V\cap G)$ is open in
$\pi_J(G)$, for every basic open neighbourhood $V$ of the identity
$e$ in $PD$. In fact, we will show that $\pi_J(V\cap G)=\pi_J(V)\cap
\pi_J(G)$.

Given a basic open neighbourhood $V$ of $e$ in $PD$, one can find
a countable set $C\subset I$ and open subgroups $V_i$ of $D_i$ for
$i\in C$ such that
\[
V=\{x\in D: x(i)\in V_i\ \mbox{for\ each}\ i\in C\}.
\]
Let $F=C\cap J$ and $E=C\setminus J$. It is clear that $F$ and $E$
are disjoint countable sets and $C=F\cup E$. Take an arbitrary point
$y\in \pi_J(V)\cap \pi_J(G)$. There exists an element $x\in G$ with
$\pi_J(x)=y$. Clearly, $x(i)=y(i)\in V_i$ for each $i\in F$. Since $E$ is
countable, we can find an element $x_0\in\Sigma{D}$ such that
$\supp(x_0)\cap J=\emptyset$ and $x_0(i)=x(i)$ for each $i\in E$.
Then the element $z=x\cdot x_0^{-1}$ of $D$ satisfies $z(i)\in V_i$
for each $i\in C$, so $z\in V$. Since $x_0\in \Sigma{D}\subset G$, we
see that $z\in G$. Hence $z\in V\cap G$ and $\pi_J(z)=\pi_J(x)\cdot
(\pi_J(x_0))^{-1}=y$. This implies that $\pi_J(V)\cap\pi_J(G)\subset
\pi_J(V\cap G)$. The inverse inclusion is obvious, so the equality
$\pi_J(V\cap G)=\pi_J(V)\cap\pi_J(G)$ is proved. Therefore the
restriction of the homomorphism $\pi_J$ to $G$ is open when
considered as a mapping onto its image.
\end{proof}

\begin{corollary}\label{Cor:N}
Let $G$ be a group with $\Sigma{D}\subset G\subset PD$, and
$\chi\in G^\wedge$. Then, for every set $J\subset I$, the continuous
extension $\bar{\chi}$ of $\chi$ over $PD$ does not depend on
$I\setminus J$ if and only if $\chi$ does not depend on $I\setminus J$.
\end{corollary}

\begin{proof}
By Lemma~\ref{Le:Dep3}, $\bar\chi$ is a continuous character on $PD$.
Hence it suffices to verify that if $\chi$ does not depend on $I\setminus J$,
neither does $\bar\chi$. Under this assumption, there exists a character
$\psi$ on $\pi_J(G)$ such that $\chi=\psi\circ\pi_J\hskip-4pt\restriction_G$.
Since the restriction to $G$ of the projection $\pi_J$ is open by
Lemma~\ref{Le:Open}, the character $\psi$ is continuous. We apply
Lemma~\ref{Le:Dep3} once again to conclude that $\bar\chi$ does not
depend on $I\setminus J$.
\end{proof}

%%%%%%%%%%%%%%%%%%%%%%%%%%%%%%%%%%
\section{Reflexivity of $G$ and compact subsets of $G^\wedge$}%%
\label{Sec:BR}%%%%%%%%%%%%%%%%%%%%%%%%%%%
Here we characterize the reflexivity of a $P$-group $G$ in terms of
compact subsets of the dual group $G^\wedge$. First we present a
simple but useful piece of information.

\begin{lemma}\label{omeb}
If $G$ is a $P$-group, then the dual group $G^\wedge$ is
$\omega$-bounded.
\end{lemma}

\begin{proof}
Since $G$ is a $P$-group, all compact subsets of $G$ are finite.
Hence $G^\wedge$ is a topological subgroup of $\T^G$, where $\T$ is
the circle group with its usual compact topology inherited from the
complex plane. Let $C$ be a countable subset of $G^\wedge$. Since
the circle group $\T$ and its power $\T^G$ are compact, the set
$\overline{C}$ is also compact (the closure is taken in $\T^G$). We
claim that $\overline{C}\subset G^\wedge$.

To verify this inclusion, take an arbitrary element $\varphi\in\overline{C}$.
It is well known (and easy to see) that $\varphi$ is a homomorphism of
$G$ to $\T$, so it suffices to prove the continuity of $\varphi$ at
the neutral element $0_G$ of $G$. Since each $\chi\in C$ is
continuous at $0_G$, there exists an open neighbourhood $U_\chi$ of
$0_G$ in $G$ such that $\chi(U_\chi)=\{1\}$; here $1$ is the neutral
element of $\T$. Then $V=\bigcap_{\chi\in C}U_\chi$ is an open
neighbourhood of $0_G$, and we claim that $\varphi(V)=\{1\}$.
Indeed, if $x\in V$, then $\chi(x)=1$ for each $\chi\in C$, so it
follows from $\varphi\in\overline{C}$ that $\varphi(x)=1$. Thus,
$\varphi$ is also continuous at $0_G$, i.e., $\varphi\in G^\wedge$.
Therefore the group $G^\wedge$ is $\omega$-bounded.
\end{proof}

Let us isolate a property that a $P$-group must possess in order
to be reflexive. In what follows we say that a set $K\subset G^\wedge$
is \textit{constant on a subgroup} $H$ of $G$ if every $\chi\in K$ is
constant on $H$.

\begin{lemma}\label{L:evm}
Let $G$ be a $P$-group. The evaluation mapping
$\alpha_{_{G}}\colon G \to G^{\wedge \wedge}$ is continuous if and only if
every compact set $K\subseteq G^\wedge$ is constant on an open
subgroup of $G$.
\end{lemma}

\begin{proof}
\emph{Necessity}. Let $K$ be a compact subset of $G^\wedge$.
If $\alpha_{_{G}}$ is continuous, there exists an open neighbourhood
$U$ of the neutral element $e$ in $G$ such that $\alpha_{_{G}}(U)
\subset K^\vartriangleright$. Since $G$ is a $P$-group, it follows from
\cite[Lemma~4.4.1\,a)]{AT} that there exists an open subgroup $V$
of $G$ such that $V\subset U$. Clearly the set $\T_+$ does not contain
nontrivial subgroups, so $\chi(V)=\{1\}$ for each $\chi\in K$.
Thus $K$ is constant on $V$.

\emph{Sufficiency}. It suffices to verify the continuity of the
homomorphism $\alpha_{_{G}}$ at the neutral element of $G$. Let
$K^\vartriangleright$ be a basic open set in $G^{\wedge\wedge}$,
with $K$ a compact subset of $G^\wedge$. By hypothesis, there
exists an open subgroup $V$ of $G$ such that every $\chi\in K$
is constant on $V$. Then $\alpha_{_{G}}(V)\subset K^\vartriangleright$.
Therefore $\alpha_{_{G}}$ is continuous at $e$.
\end{proof}

Since we are mainly concerned with (subgroups of) product groups,
it is worth to reformulate the above lemma for this special case.

\begin{corollary}\label{evm}
Let $G$ be a subgroup of $PD$, where $D=\prod_{i\in I}D_i$ is the
product of an arbitrary family of topological groups. The evaluation
mapping $\alpha_{_{G}}\colon G \to G^{\wedge \wedge}$ is continuous
if and only if for every compact set $K\subseteq G^\wedge$, one can
find a set $J\in [I]^{\leq\omega}$ and an open subgroup $U$ of
$PD_J$ such that every $\chi\in K$ is constant on the set
$G\cap\pi_J^{-1}(U)$.
\end{corollary}

\begin{proof}
By Lemma~\ref{L:evm}, the continuity of $\alpha_{_{G}}$ means that
every compact set $K\subset G^\wedge$ is constant on an open subgroup
$V$ of $G$. Since the sets $G\cap \pi_J^{-1}(U)$, where $J\in [I]^{\leq\omega}$
and $U$ is an open subgroup of $PD_J$, form an open basis at the neutral
element of $G$, the required conclusion is immediate.
\end{proof}

\begin{theorem}\label{T:duality}
A $P$-group $G$ is reflexive if and only if every compact set
$K\subset G^\wedge$ is constant on an open subgroup of $G$.
\end{theorem}

\begin{proof}
By \cite[Lemma~4.4.1\,a)]{AT}, every $P$-group has a base at the
identity consisting of open subgroups. Hence $G$ is a topological
subgroup of a product of discrete groups. Since the class of nuclear
groups contains discrete Abelian groups and is closed under taking
products and arbitrary subgroups (here we apply Propositions~7.3,
7.5, and~7.6 from \cite{bana91}), it follows that the group $G$ is nuclear.

We already know that the evaluation homomorphism $\alpha_{_{G}}
\colon G\rightarrow G^{\wedge \wedge}$ is injective because $G$
is a nuclear group. By the same reason, the mapping $\alpha_{_{G}}
\colon G\to\alpha_{_{G}}(G)$ is open, where $\alpha_{_{G}}(G)$ carries
the topology inherited from $G^{\wedge\wedge}$ (see Theorem~8.5 and
Lemma~14.3 of \cite{bana91}).

Since the $P$-group $G$ has no infinite compact subsets, $G^\wedge$
carries the topology of pointwise convergence on elements of $G$. It
follows that $G^{\wedge\wedge}=\alpha_{_{G}}(G)$ (see for instance
\cite[Theorem~1.3]{comfross64}). The $P$-group $G$ is therefore reflexive
if and only if $\alpha_{_{G}}$ is continuous. The theorem is then a direct
consequence of Lemma~\ref{L:evm}.
\end{proof}

Here is a coordinatewise form of Theorem~\ref{T:duality} which
is immediate after Corollary~\ref{evm}.

\begin{theorem}\label{duality}
Let $D=\prod_{i\in I}D_i$ be a product of topological groups. A subgroup
$G$ of $PD$ is reflexive if and only if for every compact set
$K\subset G^\wedge$, there exist a set $J\in [I]^{\leq\omega}$ and
an open subgroup $U$ of $PD_J$ such that $K$ is constant on
$G\cap\pi_J^{-1}(U)$.
\end{theorem}

The following two somewhat unexpected facts fail to hold outside
the class of $P$-groups. The first of them says that reflexivity in
$P$-groups extends from a dense subgroup to the whole group:

\begin{proposition}\label{Pro:Une}
Let $G$ be a a dense subgroup of a topological group $H$.
If $G$ is a reflexive $P$-group, so is $H$.
\end{proposition}

\begin{proof}
Suppose that $G$ is a reflexive $P$-group. Since $G$ is
dense in $H$, it follows from \cite[Lemma~4.4.1\,d)]{AT} that $H$ is a
$P$-group. According to Theorem~\ref{T:duality} it suffices to verify
that every compact set $K\subset H^\wedge$ is constant on an open
subgroup of $H$. Denote by $r$ the natural restriction homomorphism
of $H^\wedge$ to $G^\wedge$ defined by $r(\chi)=\chi\hskip-3pt\restriction_G$,
for each $\chi\in H^\wedge$. Since $G\subset H$ and the dual groups
$G^\wedge$ and $H^\wedge$ carry the topologies of pointwise
convergence on elements of $G$ and $H$, respectively, $r$ is continuous.
It also follows from the density of $G$ in $H$ that $r$ is one-to-one and onto.
In other words, $r$ is a continuous isomorphism of $H^\wedge$ onto $G^\wedge$.

Since $G$ is reflexive, Theorem~\ref{T:duality} implies that the
compact set $r(K)\subset G^\wedge$ is constant on an open subgroup
$V$ of $G$. Let $U$ be the closure of $V$ in $H$. Then $U$ is an
open subgroup of $H$ and every $\chi\in K$ is constant on the dense
subset $V$ of $U$. By a continuity argument, $\chi$ is constant on
$U$. Therefore $K$ is constant on $U$, whence the reflexivity of $H$
follows.
\end{proof}

The second fact establishes that the class of reflexive $P$-groups
is stable under taking quotients. Its proof makes use of \textit{dual
homomorphisms.} Since this tool will be used several times in the
article, we give a lemma explaining basic properties of dual
homomorphisms.

Let us recall that a surjective mapping $f\colon X\to Y$ is
\textit{compact covering\/} if for every compact set $K\subset Y$,
there exists a compact set $C\subset X$ such that $f(C)=K$.

\begin{lemma}\label{Le:Closed}
Let $\pi\colon G\to H$ be a continuous homomorphism of
topological Abelian groups. Let also $\pi^\wedge\colon
H^\wedge\to G^\wedge$ be the dual homomorphism defined
by $\pi^\wedge(\chi)=\chi\circ\pi$, for each $\chi\in H^\wedge$.
Then:
\begin{enumerate}
\item[\rm (a)] $\pi^\wedge$ is continuous.
\item[\rm (b)] If $\pi$ is compact covering, then $\pi^\wedge$
is a topological isomorphism of $H^\wedge$ onto its image
$\pi^\wedge(H^\wedge)$.
\item[\rm (c)] If the homomorphism $\pi$ is open, then the image
$\pi^\wedge(H^\wedge)$ is closed in $G^\wedge$.
\end{enumerate}
\end{lemma}

\begin{proof}
(a) is well known. Indeed, let $C$ be a compact subset of $G$
and $U=C^\vartriangleright$ a basic open neighbourhood of the
neutral element in $G^\wedge$. Then $K=\pi(C)$ is a compact
subset of $H$ and $V=K^\vartriangleright$ is an open
neighbourhood of the neutral element in $H^\wedge$. For any
$\chi\in V$, we have $\pi^\wedge(\chi)(C)=(\chi\circ\pi)(C)=
\chi(\pi(C))=\chi(K)\subset \T_+$, that is, $\pi^\wedge(\chi)\in C
^\vartriangleright=U$. This implies the continuity of $\pi^\wedge$ and
proves (a).\smallskip

(b) follows from \cite[Lemma~5.17]{Aub}.\smallskip

(c) Suppose that $\pi$ is an open homomorphism of $G$ to $H$
and let $K=\pi(G)$. Then $K$ is an open subgroup of $H$, so every
continuous character on $K$ extends to a continuous character on
$H$. Hence $\pi^\wedge(K^\wedge)=\pi^\wedge(H^\wedge)$ and we
can assume without loss of generality that $\pi(G)=H$.

Our further argument is very close to the proof of item 2) of Corollary~0.4.8
from \cite{Arh}. Indeed, let $\psi$ be in the closure of $\pi^\wedge
(H^\wedge)$ in $G^\wedge$. Since finite sets are compact, the
topology of $G^\wedge$ contains the topology of pointwise
convergence on elements of $G$. Applying this fact, one easily
verifies that $\psi$ is a homomorphism of $G$ to $\T$ and that $\psi$
is constant on each fiber $\pi^{-1}(y)$, $y\in H$. Hence there exists
a function $\chi\colon H\to\T$ such that $\psi=\chi\circ\pi$. Clearly,
$\chi$ is a homomorphism. Since $\psi$ is continuous and $\pi$ is
open and onto, the equality $\psi=\chi\circ\pi$ implies that $\chi$ is
continuous as well. Thus, $\chi\in H^\wedge$ and
$\psi=\pi^\wedge(\chi)\in\pi^\wedge(H^\wedge)$.
\end{proof}

\begin{corollary}\label{Fins}
Let $\pi\colon G\to H$ be a continuous onto homomorphism.
If all compact subsets of $H$ are finite, then $\pi^\wedge$ is a
topological isomorphism of $H^\wedge$ onto the subgroup
$\pi^\wedge(H^\wedge)$ of $G^\wedge$. In particular, this is
the case when $H$ is a $P$-group.
\end{corollary}

\begin{proof}
Since all compact subsets of $H$ are finite, the homomorphism
$\pi$ is compact covering. The required conclusion now follows
from (b) of Lemma~\ref{Le:Closed}.
\end{proof}

\begin{proposition}\label{Pro:Quo}
Let $\pi\colon G\to H$ be a continuous open epimorphism of topological
groups. If $G$ is a reflexive $P$-group, so is $H$.
\end{proposition}

\begin{proof}
The fact that the image $H=\pi(G)$ is a $P$-group follows from
\cite[Lemma~4.4.1\,c)]{AT}. Let us show that $H$ is reflexive. Take an
arbitrary compact subset $C$ of $H^\wedge$. By (a) of Lemma~\ref{Le:Closed},
the dual homomorphism $\pi^\wedge\colon H^\wedge\to G^\wedge$ is
continuous, so $K=\pi^\wedge(C)$ is a compact subset of $G^\wedge$.
Since the group $G$ is reflexive, Theorem~\ref{T:duality} implies that
there exists an open subgroup $U$ of $G$ such that $K$ is constant
on $U$. Then $C$ is constant on the open subgroup $V=\pi(U)$ of $H$.
Indeed, if $\chi\in C$ and $y\in V$, take $x\in U$ with $\pi(x)=y$. Then
$\chi(y)=\chi(\pi(x))=\pi^\wedge(\chi)(x)=1$ since $\pi^\wedge(\chi)\in K$.
Applying Theorem~\ref{T:duality} once again, we conclude that the group
$H$ is reflexive.
\end{proof}

We finish this section with two special cases of Theorem~\ref{duality}.

\begin{corollary}\label{CorSp1}
Let $D=\prod_{i\in I}D_i$ be a product of topological groups. The group
$PD$ is reflexive if and only if the following hold:
\begin{enumerate}
\item[{\rm (a)}] the group $PD_J=P(\prod_{i\in J}D_i)$ is reflexive for
each $J\in [I]^{\leq\omega};$
\item[{\rm (b)}] every compact set $K\subset (PD)^\wedge$ depends
at most on countably many coordinates.
\end{enumerate}
\end{corollary}

\begin{proof}
\textit{Necessity.} It is easy to see that the projection $\pi_J\colon PD\to PD_J$
is open for each $J\subset I$. If the group $PD$ is reflexive, then the reflexivity
of the groups $PD_J$ follows from Proposition~\ref{Pro:Quo}. Suppose that $K$
is a compact subset of $(PD)^\wedge$. Theorem~\ref{duality} implies that there
exist $J\in [I]^{\leq\omega}$ and an open subgroup $U$ of $PD_J$ such that
$K$ is constant on $\pi_J^{-1}(U)$. Then every $\chi\in K$ does not depend
on $I\setminus J$, that is, $K$ does not depend on $I\setminus J$. Hence
conditions (a) and (b) hold true.

\textit{Sufficiency.} Let us deduce the reflexivity of $PD$ from (a) and (b).
Take a compact set $K\subset (PD)^\wedge$. By (b), there exists a countable
set $J\subset I$ such that $K$ does not depend on $I\setminus J$.

Denote by $\varphi\colon (PD_J)^\wedge \to (PD)^\wedge$ the
homomorphism dual to the projection $\pi_J\colon PD\to PD_J$. Since
the projection $\pi_J\colon PD\to PD_J$ is open and all compact
subsets of $PD_J$ are finite, it follows from Corollary~\ref{Fins}
that $\varphi$ is a topological isomorphism of $(PD_J)^\wedge$ onto
a subgroup of $(PD)^\wedge$.

We claim that $K\subset \varphi((PD_J)^\wedge)$. Indeed, take an arbitrary
character $\chi\in K$. Since $\chi$ does not depend on $I\setminus J$, there
exists a character $\zeta$ on $PD_J$ such that $\chi=\zeta\circ\pi_J$. The
character $\zeta$ is continuous since the projection $\pi_J$ is open. Hence
$\zeta\in (PD_J)^\wedge$ and $\chi=\varphi(\zeta)\in \varphi((PD_J)^\wedge)$.

Let $C=\varphi^{-1}(K)$. Then $C$ is a compact subset of the group
$(PD_J)^\wedge$. By (a), the group $PD_J$ is reflexive. According to
Theorem~\ref{T:duality}, $PD_J$ contains an open subgroup $U$ such
that every character $\zeta\in C$ is constant on $U$. Since $\varphi(C)=K$,
we see that every $\chi\in K$ is constant on the open subgroup $\pi_J^{-1}(U)$
of $PD$. The reflexivity of the group $PD$ now follows from Theorem~\ref{duality}.
\end{proof}

We shall see in Proposition~\ref{Pro:Ge} that one can drop item (b) in the above
corollary. The next result is a slight modification of Corollary~\ref{CorSp1}, so
we omit its proof.

\begin{corollary}\label{CorSp2}
Let $D=\prod_{i\in I}D_i$ be a product of topological groups such
that the group $PD_J$ is reflexive for each $J\in [I]^{\leq\omega}$.
Suppose that $G$ is a subgroup of $PD$ satisfying $\pi_J(G)=PD_J$
for each $J\in [I]^{\leq\omega}$. Then $G$ is reflexive if and only if
every compact set $K\subset G^\wedge$ depends at most on countably
many coordinates.
\end{corollary}

%%%%%%%%%%%%%%%%%%%%%%%%%%%%%%%%%%%%%
\section{Reflexive $P$-groups}\label{RPG}%%%%%%%%%%%%%%%%
We prepare here our way to show that some nondiscrete $P$-groups are
reflexive. The lemma below is obvious and its proof is omitted.

\begin{lemma}\label{folk}
Let $b,t$ be elements of $\T$ and $t\neq 1$. Then there is an integer
$k$ such that $t^k\cdot b\notin \T_{+}$.
\end{lemma}

In the following proposition and in Lemmas~\ref{triv}--\ref{Le:MT},
$D=\prod_{i\in I} D_i$ stands for the product of an arbitrary family of topological
Abelian groups.

\begin{proposition}\label{construction}
Suppose that $C=\left\{ \chi_\eta\colon \eta< \omega_1\right\}\subset
(PD)^\wedge$, $\mathcal{J}=\{J_\eta \colon\eta<\omega_1\}\subset
[I]^{\leq \omega}$, and $X=\{x_\eta \colon \eta<\omega_1\}\subset PD$
satisfy the following conditions for each $\eta< \omega_1$:
\begin{enumerate}
\item $\chi_\eta$ does not depend on $I\setminus J_\eta;$
\item $\supp(x_\eta)\subset J_\eta ;$
\item $\chi_\eta(x_\eta)\in\T\setminus\T_{+};$
\item if $\,\zeta<\eta$, then $J_\zeta\bigcap \supp(x_\eta)=
\emptyset$.
\end{enumerate}
Then every element of $\,\bigcap_{\gamma<\omega_1} \overline{\left\{
\chi_\eta \colon \eta\geq \gamma\right\}}^{\,\T^D}$ is discontinuous
as a character on the group $PD$.
\end{proposition}

\begin{proof}
For each $\gamma<\omega_1$, let $C_\gamma=\{\chi_\eta \colon \eta
\geq \gamma\}$ and $K_\gamma=\overline{C_\gamma}^{\hskip 2pt\T^D}$.
Since the family $\{K_\gamma \colon \gamma < \omega_1\}$ of compact
sets is decreasing, we see that $K=\bigcap_{\gamma <\omega_1} K_\gamma$
is nonempty. Take $\rho\in K$. Suppose toward a contradiction that $\rho$
is continuous, i.e., that $\rho\in (PD)^\wedge$. By Lemma~\ref{Le:Depen},
there exists $J\in [I]^{\leq\omega}$ such that $\rho$ does not depend on
$I\setminus J$. It follows from (2) and (4) that the family
$\{\supp(x_\eta)\colon \eta <\omega_1\}$ is pairwise disjoint.
Take a countable ordinal $\eta_0$ such that
$J \bigcap\supp(x_\eta)=\emptyset$ for all $\eta$ satisfying
$\eta_0\leq \eta<\omega_1$. It then follows that $\rho (x_\eta)=1$
for every countable ordinal $\eta \geq \eta_0$.

Given a family $\{g_\alpha: \alpha\in A\}\subset D$ such that
$\supp(g_\alpha)\cap\supp(g_\beta)=\emptyset$ for distinct
$\alpha,\beta\in A$, we can define an element $g=\coprod_{\alpha\in
A}g_\alpha\in D$ by the requirements that $\supp(g)=
\bigcup_{\alpha\in A}\supp(g_\alpha)$ and for every $\alpha\in A$,
the elements $g$ and $g_\alpha$ coincide on $\supp(g_\alpha)$.

Now for every countable ordinal $\eta \geq \eta_0$, we define a
point $g_\eta\in D$ satisfying the following two conditions:
\begin{enumerate}
\item[(a)] $g_\eta$ is either the neutral element $\mathbf 0$ of
               $D$ or $k_\eta x_\eta$, for some $k_\eta \in \Z$;
\item[(b)] $\chi_\eta\left(\coprod_{\eta_0\leq \beta\leq \eta}
               g_\beta\right)\in\T\setminus \T_{+}$.
\end{enumerate}
To begin, we put $g_{\eta_0}=x_{\eta_0}$. Then conditions (a) and
(b) hold. Suppose now that $\eta_0 < \sigma < \omega_1$ and that
$g_\eta$ have been defined for all $\eta$ with $\eta_0 < \eta<\sigma$
such that (a) and (b) hold. Notice that the family
$\{\supp(g_\eta): \eta_0\leq\eta<\sigma\}$ is pairwise disjoint, by
(a). If $\chi_\sigma\left(\coprod_{\eta_0\leq\beta< \sigma}
g_\beta \right)\notin\T_{+}$, put $g_\sigma=\mathbf{0}$. If
$z=\chi_\sigma\left(\coprod_{\eta_0\leq\beta< \sigma}
g_\beta \right)\in \T_+$, we apply Lemma~\ref{folk} to
$z$ and $t=\chi_\sigma (x_\sigma)$ to find an integer $k_\sigma$
such that $z\cdot t^{k_\sigma}\notin\T_{+}$. We then put
$g_\sigma =k_\sigma x_\sigma$. Since
$\chi_\sigma\left(\coprod_{\eta_0\leq\beta\leq\sigma}g_\beta\right)=
z\cdot t^{k_\sigma}\notin\T_{+}$, the element $g_\sigma$
satisfies (a) and (b) at the stage $\sigma$. The recursive
definitions are complete.

By (a), the supports of $g_{\eta}$'s with $\eta\geq\eta_0$ are
disjoint, so we can put $h_0=\coprod_{\eta_0\leq \eta<\omega_1}
g_\eta$. Again by (a), $\supp(h_0)\subset \bigcup
\{\supp(x_\eta)\colon \eta_0 \leq \eta < \omega_1\}$ and, since
$J\cap \supp(x_\eta)= \emptyset$ for every countable ordinal $\eta
\geq \eta_0$, we see that $J\cap \supp(h_0)=\emptyset$. Therefore
$\rho(h_0)=1$. Since $\rho \in K_{\eta_0}$ and $(\T_+)_{h_0}\times
\T^{D\setminus\{h_0\}}$ is a basic open set in $\T^D$ containing
$\rho$, there exists $\eta \geq \eta_{0}$ such that
$\chi_{\eta}(h_0)\in \T_+$. However, it follows from (4) of the
proposition that $J_\eta\cap\supp(x_\beta)=\emptyset$ if
$\eta<\beta<\omega_1$, while condition (a) implies that
$\supp(g_\beta)\subset\supp(x_\beta)$. Therefore, the sets
$J_\eta$ and $\supp\left( \coprod_{\eta<\beta<\omega_1} g_\beta\right)$
are disjoint and hence (1) of the proposition implies that
$\chi_\eta\left(\coprod_{\eta<\beta<\omega_1} g_\beta\right)=1$.
It now follows from (b) that
\[
\chi_\eta(h_0)=\chi_\eta\left(\coprod_{\eta_0\leq\beta\leq \eta} g_\beta\right)
\cdot\chi_\eta\left(\coprod_{\eta<\beta<\omega_1}g_\beta\right)=
\chi_\eta\left(\coprod_{\eta_0\leq\beta\leq \eta }g_\beta\right)\cdot
1\notin\T_{+}.
\]
This contradiction shows that $\rho$ is discontinuous.
\end{proof}

\begin{lemma}\label{triv}
Suppose that $\Sigma{D}\subset G\subset PD$, $g\in G$, $\chi\in
G^\wedge$, and $\chi(g)\neq 1$. If $\chi$ does not depend on
$I\setminus J$, for a countable set $J\subset I$, then there exists
a point $x\in G$ with $\supp(x)\subset J\cap\supp(g)$ such that
$\chi(x)\notin\T_{+}$.
\end{lemma}

\begin{proof}
We can find elements $y,z\in PD$ such that $g=y+z$, $\supp(y)=
J\cap\supp(g)$, and $\supp(z)\cap J=\emptyset$. Notice that
$y\in\Sigma{D}\subset G$, so $z\in G$. It follows from our choice of
$z$ that $\chi(z)=1$. Since $\chi(y)=\chi (y)\cdot \chi(z)=\chi(y+z)=
\chi(g)\neq 1$, we see that $t=\chi(y)\neq 1$. Take an integer $n$
such that $t^n\notin\T_{+}$. The point $x=ny\in G$ is as required.
\end{proof}

\begin{lemma}\label{main}
Suppose that $\Sigma{D}\subseteq G \subseteq PD$ and $K\subseteq
G^\wedge$. If $K$ depends on uncountably many coordinates, then $K$
contains a subset $C=\{ \chi_\eta\colon \eta<\omega_1\}$ and $G$
contains a subset $X=\{x_\eta \colon \eta<\omega_1\}$ such that
conditions (1)--(4) of Proposition~\ref{construction} hold for $C$, $X$,
and a suitable collection $\mathcal{J}=\{J_\eta\colon \eta < \omega_1 \}$
of countable subsets of $I$.
\end{lemma}

\begin{proof}
To begin, we choose a family $\{J_\chi\colon \chi \in K\}$ of countable
subsets of $I$ such that $\chi$ does not depend on $I\setminus J_\chi$,
for each $\chi\in K$.

Since $K$ depends on uncountably many coordinates, there is a
nontrivial character $\chi_0 \in K$. It follows from Lemma~\ref{triv}
that there is a point $x_0\in G$ with $\supp(x_0)\subset J_{\chi_0}$
such that $\chi_0(x_0)\notin\T_{+}$. We put $J_0=J_{\chi_0}$.
Then (1)--(4) of Proposition~\ref{construction} hold for $\eta=0$.

Suppose now that $\sigma < \omega_1$ and that $\chi_\eta$,
$x_\eta$, and $J_\eta$ have been defined to satisfy (1)--(4) of
Proposition~\ref{construction} for all $\eta <\sigma$. Then the set
$T=\bigcup_{\eta <\sigma}J_\eta$ is countable. By assumptions of
the lemma, there is a character $\chi_\sigma \in K$ depending on
$I\setminus T$, i.e., there is a point $g\in G$ such that $g(i)$ is the
neutral element of $D_i$ for each $i\in T$ and $\chi_\sigma (g)\neq 1$.
It now follows from Lemma~\ref{triv} that there exists a point
$x_\sigma\in G$ with $\supp(x_\sigma)\subset \supp(g)\cap J_{\chi_\sigma}$
such that $\chi_\sigma(x_\sigma)\notin\T_{+}$. Let $J_\sigma=J_{\chi_\sigma}$.
It is clear that $\chi_\eta$, $x_\eta$, and $J_\eta$ satisfy (1)--(4) of
Proposition~\ref{construction} for all $\eta\leq\sigma$. The recursive
definitions are complete which finishes the proof.
\end{proof}

A slight modification in the above argument can be made in order
to deduce the following lemma which will be applied in the proof
of Proposition~\ref{main1}.

\begin{lemma}\label{Le:MT}
Suppose that $\Sigma{D}\subset G\subset PD$ and that two sets $C=\{\chi_\eta:
\eta<\omega_1\}\subset G^\wedge$ and $\mathcal{J}=\{J_\eta: \eta<\omega_1\}
\subset [I]^{\leq\omega}$ satisfy the following conditions for all $\eta<\omega_1$:
\begin{enumerate}
\item[\rm (a)] $\chi_\eta$ does not depend on $I\setminus J_\eta$;
\item[\rm (b)] $J_\zeta\subset J_\eta$ if $\zeta<\eta$;
\item[\rm (c)] $\chi_\eta$ depends on $I\setminus \bigcup_{\zeta<\eta}J_\zeta$.
\end{enumerate}
Then there exists a set $X=\{x_\eta: \eta<\omega_1\}\subset G$ such that
$C$, $\mathcal{J}$, and $X$ satisfy conditions (1)--(4) of
Proposition~\ref{construction}.
\end{lemma}

\begin{proposition}\label{Pro:Count}
Let $D=\prod_{i\in I}D_i$ be a product of topological groups.
Then every compact set $K\subset (PD)^\wedge$ depends
at most on countably many coordinates.
\end{proposition}

\begin{proof}
Let $\Pi=PD$. Suppose to the contrary that a compact set
$K\subset\Pi^\wedge$ depends on uncountably many coordinates.
By Lemma~\ref{main}, there exist $C\subset K$, $X\subset G$,
and $\mathcal{J}\subset [I]^{\leq\omega}$ satisfying conditions
(1)--(4) of Proposition~\ref{construction}. By the latter proposition,
there is an element $\rho$ in the closure of $C$ in $\T^{\Pi}$ such that
$\rho\notin \Pi^\wedge$. But $K$, being compact, must be closed
in $\T^\Pi$, whence $\rho\in K\subset\Pi^\wedge$. This contradiction
shows that there must exist $J\in [I]^{\leq\omega}$ such that $K$ does
not depend on $I\setminus J$.
\end{proof}

The above proposition shows that item (b) in Corollary~\ref{CorSp1}
can be omitted:

\begin{proposition}\label{Pro:Ge}
Let $D=\prod_{i\in I}D_i$ be a product of topological groups.
Then the product group $PD$ is reflexive if and only if $PD_J$
is reflexive for each $J\in [I]^{\leq\omega}$.
\end{proposition}

\begin{theorem}\label{cor:main}
Let $D=\prod_{i\in I}D_i$ be a product of discrete Abelian groups.
Then the $P$-group $\Pi=PD$ and the $\omega$-bounded group
$\Pi^\wedge$ are reflexive.
\end{theorem}

\begin{proof}
The group $PD_J$ is discrete and hence reflexive, for every
$J\in [I]^{\leq\omega}$. Therefore the reflexivity of $\Pi$ follows
from Proposition~\ref{Pro:Ge}. Hence the dual group $\Pi^\wedge$
is reflexive as well.
\end{proof}

In the case when the product $D=\prod_{i\in I}D_i$ in the above theorem
contains uncountably many nontrivial factors, we obtain the following
result that answers a question posed by S.~Hern\'andez and P.~Nickolas
and solves a problem raised in a comment after Proposition~2.10 in
\cite{ardaetal}.

\begin{corollary}\label{Cor:Exis}
There exist non-discrete reflexive $P$-groups as well as $\omega$-bounded
noncompact reflexive groups.
\end{corollary}

We can now establish the reflexivity of certain $P$-groups which are not
necessarily $P$-modifications of products of discrete groups. A simple
auxiliary lemma is in order:

\begin{lemma}\label{Le:Aux}
Suppose that $\pi\colon G\to H$ is a continuous onto homomorphism of
compact groups. Then the homomorphism $\pi\colon PG\to PH$ is open,
where $PG$ and $PH$ are $P$-modifications of the groups $G$ and $H$,
respectively.
\end{lemma}

\begin{proof}
Let $e$ be the neutral element of $G$. It is clear that the sets of the form
$V=\bigcap_{n\in\omega}U_n$, where $U_n$'s are open neighbourhoods
of $e$ in $G$ and $\overline{U}_{n+1}\subset U_n$ for each $n\in\omega$
(the closure is taken in $G$), constitute a base at $e$ in $PG$. Therefore,
it suffices to verify that every image $\pi(V)$ is open in $PH$. Notice that
the continuous epimorphism $\pi\colon G\to H$ is open since $G$ is compact.
Using the compactness of $G$ once again we see that $\pi(V)=\bigcap
_{n\in\omega}\pi(U_n)$, so $\pi(V)$ is a $G_\delta$-set in $H$. Hence
$\pi(V)$ is open in $PH$.
\end{proof}

\begin{proposition}\label{Cor:Comp}
Let $H$ be a compact Abelian group and $PH$ the $P$-modification
of $H$. Then the group $PH$ is reflexive.
\end{proposition}

\begin{proof}
It is well known that one can find a compact Abelian group $G$ of the form
$G=\prod_{i\in I}G_i$, with compact metrizable factors $G_i$, and a continuous
homomorphism $\pi$ of $G$ onto $H$ (see \cite[Lemma~1.6]{HM}). By
Lemma~\ref{Le:Aux}, the homomorphism $\pi\colon PG\to PH$ is open.
For every $i\in I$, let $D_i$ be the group $G_i$ with the discrete topology.
Denote by $D$ the product group $\prod_{i\in I}D_i$. Since the factors
$G_i$ are metrizable, the topological groups $PG$ and $PD$ coincide.
It now follows from Theorem~\ref{cor:main} that the group $PG$ is reflexive,
while Proposition~\ref{Pro:Quo} implies the reflexivity of $PH$.
\end{proof}

According to \cite[Section~4.3]{AT}, a topological group $H$ is \textit{feathered\/}
if it contains a nonempty compact set with a countable neighbourhood base in
$H$. Let us call a topological group $H$ \textit{pseudo-feathered\/} if there exists
a nonempty compact set of type $G_\delta$ in $H$. It is clear that every
feathered group is pseudo-feathered and that $H$ is pseudo-feathered if and
only if it contains a compact subgroup of type $G_\delta$. An Abelian group
is pseudo-feathered iff it admits an open continuous homomorphism with
compact kernel onto a group of countable pseudocharacter. In the following
result we extend the conclusion of Proposition~\ref{Cor:Comp} to
pseudo-feathered groups.

\begin{proposition}\label{Pro:Ex}
Let $H$ be a pseudo-feathered Abelian group. Then the group $PH$
is reflexive.
\end{proposition}
%\begin{proof}
%Suppose that $C$ is a compact subgroup of type $G_\delta$ in $H$.
%Clearly, $C$ is open in $PH$. Let  $K$ be a compact subset of
%$(PH)^\wedge$. Denote by $r$ the restriction mapping of
%$(PH)^\wedge$ to $(PC)^\wedge$,
%$r(\chi)=\chi\hskip-3pt\restriction_C$ for each $\chi\in
%(PH)^\wedge$. It is clear that $r$ is continuous, so $r(K)$ is a
%compact subset of $(PC)^\wedge$. Since the group $PC$ is reflexive
%by Proposition~\ref{Cor:Comp}, we can apply Theorem~\ref{T:duality}
%to find an open subgroup $U$ of $PC$ such that every $\xi\in r(K)$
%is constant on $U$. Then every element $\chi\in K$ is constant on
%$U$, while $U$ is open in $PH$. Now the reflexivity of $PH$ $%follows
%·from Theorem~\ref{T:duality}.
%\end{proof}
\begin{proof}
  Let $C$ be a compact subgroup of type $G_\delta$ in $H$. Clearly,
$PC$ is then an  open subgroup of $PH$. Since a topological group
admitting an open subgroup that is reflexive is itself reflexive
(Proposition 2.2 of \cite{banachasmart}) and the group $PC$ is
reflexive by Proposition~\ref{Cor:Comp}, we conclude that $PH$ is
reflexive.
\end{proof}

The next result is a common generalization of Theorem~\ref{cor:main}
and Proposition~\ref{Pro:Ex}:

\begin{theorem}\label{Th:Gen}
Let $H=\prod_{i\in I} H_i$ be the product of a family of pseudo-feathered
Abelian groups. Then the group $PH$ is reflexive.
\end{theorem}

\begin{proof}
It is easy to verify that if $C_n$ is a compact set of type $G_\delta$ in a
space $X_n$, for each $n\in\omega$, then the compact set $C=\prod
_{n\in\omega} C_n$ has type $G_\delta$ in the product space $X=\prod
_{n\in\omega} X_n$. This observation implies that the group $H_J=
\prod_{i\in J} H_i$ is pseudo-feathered for each $J\in [I]^{\leq\omega}$.
The reflexivity of $PH$ now follows from Proposition~\ref{Pro:Ge}.
\end{proof}

In Theorem~\ref{Reduc} below we characterize the reflexivity of certain
subgroups $G$ of \lq\lq{big\rq\rq} products $PD=P\prod_{i\in I}D_i$ of
topological groups in terms of projections $\pi_J(G)$ of $G$ to relatively
small subproducts $PD_J=P\prod_{i\in J}D_i$.

\begin{theorem}\label{Reduc}
Suppose that $D=\prod_{i\in I}D_i$ is a product of topological groups
and $\Sigma{D}\subset G\subset PD$. Then the group $G$ is reflexive
iff the subgroup $\pi_J(G)$ of $PD_J$ is reflexive, for every set $J\subset I$
satisfying $|J|\leq\aleph_1$.
\end{theorem}

\begin{proof}
\emph{Necessity.} Let $G$ be reflexive. Take any $J\subset I$ satisfying
$|J|\leq\aleph_1$ and put $H=\pi_J(G)$. By Lemma~\ref{Le:Open},
the restriction to $G$ of the projection $\pi_J\colon PD\to PD_J$ is an
open homomorphism of $G$ onto $H$. Hence the reflexivity of $H$
follows from Proposition~\ref{Pro:Quo}.

\emph{Sufficiency.} Suppose that $\pi_J(G)$ is reflexive, for each
$J\subset I$ with $|J|\leq\aleph_1$. Since $\Sigma{D}\subset G$,
the equality $\pi_J(G)=D_J$ holds for each $J\in [I]^{\leq\omega}$.
Therefore, according to Corollary~\ref{CorSp2}, it suffices to show
that every compact set $K\subset G^\wedge$ depends at most on
countably many coordinates. Suppose to the contrary that $G^\wedge$
contains a compact set $K$ which depends on uncountably many
coordinates. Apply Lemma~\ref{main} to choose families $\{\chi_\eta:
\eta<\omega_1\}\subset K$, $\{x_\eta: \eta<\omega_1\}\subset G$,
and $\{J_\eta: \eta<\omega_1\}\subset [I]^{\leq\omega}$ satisfying
conditions (1)--(4) of Proposition~\ref{construction}. Let $J=\bigcup
_{\eta<\omega_1}J_\eta$. Then $J\subset I$ and $|J|\leq\aleph_1$.
Hence the subgroup $H=\pi_J(G)$ of $PD_J$ is reflexive. Notice that
by Lemma~\ref{Le:Open}, the restriction to $G$ of the homomorphism
$\pi_J$ is open when considered as a mapping of $G$ onto $H$.
Let $\eta<\omega_1$. Since $\chi_\eta$ does not depend on
$I\setminus J_\eta$ and $J_\eta\subset J$, there exists a continuous
character $\psi_\eta$ on $H$ such that $\chi_\eta=\psi_\eta\circ\pi_J
\hskip-4pt\restriction_G$. We put $\Psi=\{\psi_\eta: \eta<\omega_1\}$.

Denote by $\varphi$ the continuous homomorphism $(\pi_J\hskip-4pt
\restriction_G)^\wedge$ of $H^\wedge$ to $G^\wedge$. Then $\varphi
(\psi_\eta)=\chi_\eta\in K$ for each $\eta<\omega_1$, so $\varphi(\Psi)
\subset K$. Since $H$ is a $P$-group, all compact subsets of $H$ are
finite. Hence Corollary~\ref{Fins} implies that $\varphi$ is a topological
isomorphism of $H^\wedge$ onto the subgroup $\varphi(H^\wedge)$
of $G^\wedge$. Further, since the homomorphism $\pi_J\hskip-4pt
\restriction_G$ of $G$ onto $H$ is open, it follows from item (c) of
Lemma~\ref{Le:Closed} that $\varphi(H^\wedge)$ is a closed subgroup
of $G^\wedge$. Therefore, $C=K\cap\varphi(H^\wedge)$ is a compact
subset of $\varphi(H^\wedge)$ and $L=\varphi^{-1}(C)$ is a compact
subset of $H^\wedge$. It follows from $\Psi\subset H^\wedge$ and
$\varphi(\Psi)\subset K$ that $\Psi\subset L$. The latter inclusion
and the definition of the set $\Psi$ together imply that $L$ depends
on uncountably many coordinates.

Indeed, suppose that for some countable set $A\subset J$, every
element of $L$ does not depend on $J\setminus A$. In particular,
$\psi_\eta$ does not depend on $J\setminus A$, for each
$\eta<\omega_1$. Since $\chi_\eta=\psi_\eta\circ\pi_J\hskip-4pt
\restriction_G$, we see that each $\chi_\eta$ does not depend
on $I\setminus A$. It follows from our choice of the families
$\{\chi_\eta: \eta<\omega_1\}$, $\{x_\eta: \eta<\omega_1\}$,
and $\{J_\eta: \eta<\omega_1\}$ (see conditions  (2)--(4) of
Proposition~\ref{construction}) that $\supp(x_\eta)\subset
J_\eta\setminus\bigcup_{\zeta<\eta}J_\zeta$ and $\chi_\eta(x_\eta)
\neq 1$,  for each $\eta<\omega_1$. Since the sets $A_\eta=
J_\eta\setminus\bigcup_{\zeta<\eta}J_\zeta$ are pairwise disjoint,
there exists $\eta<\omega_1$ such that $A\cap A_\eta=\emptyset$.
Since $\chi_\eta$ does not depend on $I\setminus A$, this implies
that $\chi_\eta(x_\eta)=1$, which is a contradiction. We have thus
proved that every compact subset $K$ of $G^\wedge$ depends at
most on countably many coordinates and, therefore, $G$ is reflexive.
\end{proof}

Theorem~\ref{Reduc} makes it possible to find many proper dense
reflexive subgroups of big products of pseudo-feathered groups
endowed with the $P$-modified topology:

\begin{corollary}\label{Cor:Sigma}
Suppose that $D=\prod_{i\in I} D_i$ is a product of pseudo-feathered
Abelian groups and let
\[
\Sigma_{\aleph_1}{D}=\{x\in PD: |\supp(x)|\leq\aleph_1\}.
\]
Then every group $G$ with $\Sigma_{\aleph_1}{D}\subset G\subset PD$
is reflexive.
\end{corollary}

\begin{proof}
According to Proposition~\ref{Pro:Une} it suffices to show that
$\Sigma_{\aleph_1}{D}$ is reflexive. It is clear that
$\pi_J(\Sigma_{\aleph_1}{D})=D_J=\prod_{i\in J}D_i$ for each
$J\subset I$ with $|J|\leq\aleph_1$, where $\pi_J\colon D\to D_J$
is the projection. By Theorem~\ref{Th:Gen}, the groups $PD_J$
are reflexive. One applies Theorem~\ref{Reduc} to conclude that
the group $\Sigma_{\aleph_1}{D}$ is reflexive as well.
\end{proof}

In what follows we identify the additive group $\Z_2=\{0,1\}$ with
the multiplicative subgroup $\{1,-1\}$ of $\T$. Hence the dual group
$G^\wedge$ of every boolean $P$-group $G$ is topologically isomorphic
to a subgroup of $\Z_2^G$. We will now show that the $P$-group
$\Pi=P\Z_2^{\omega_1}$ contains proper dense reflexive subgroups
of the form $G_\xi$ defined after Theorem~\ref{sigmal}.

\begin{proposition}\label{main1}
Every compact set $K\subset (G_\xi)^\wedge$ depends at most on
countably many coordinates, where $\xi$ is an ultrafilter on $\omega_1$
containing the complements to countable sets. Hence the group $G_\xi$
is reflexive.
\end{proposition}

\begin{proof}
On the contrary, suppose that a compact set $K\subset (G_\xi)^\wedge$
depends on uncountably many coordinates. We construct two sets
$\{\chi_\eta: \eta<\omega_1\}\subset K$ and $\{J_\eta \colon
\eta<\omega_1\}\subset [\omega_1]^{\leq \omega}$  satisfying the
following conditions for all $\eta<\omega_1$:
\begin{enumerate}
\item[(i)]    $\chi_\eta$ does not depend on $\omega_1\setminus J_\eta$;
\item[(ii)]   $J_\zeta\subset J_\eta$ if $\zeta<\eta$;
\item[(iii)]  $\eta\in J_\eta$;
\item[(iv)]  $\chi_\eta$  depends on the set $\omega_1\setminus
                 \left(\{\eta\}\cup\bigcup_{\zeta<\eta}J_\zeta\right)$.
\end{enumerate}
Let $\chi_0\in K$ be a nontrivial character. Take a countable set
$J_0\subset \omega_1$ such that $0\in J_0$ and $\chi_0$ does not
depend on $\omega_1\setminus J_0$. Suppose that for some
$\eta<\omega_1$, the sequences $\{\chi_\zeta:\zeta<\eta\}\subset K$
and $\{J_\zeta: \zeta<\eta\}\subset[\omega_1]^{\leq\omega}$ have
been defined to satisfy conditions (i)--(iv). Then we put $T_\eta=
\bigcup_{\zeta<\eta}J_\zeta$ and choose $\chi_\eta\in K$ such that
$\chi_\eta$ depends on the set
$\omega_1\setminus\left(T_\eta\cup\{\eta\}\right)$. Such a choice of
$\chi_\eta$ is possible since the set $T_\eta\cup \{\eta\}$ is
countable and $K$ depends on uncountably many coordinates. Let
$J'_\eta$ be a countable subset of $\omega_1$ such that $\chi_\eta$
does not depend on $\omega_1\setminus J'_\eta$. Then the set
$J_\eta=J'_\eta\cup T_\eta\cup\{\eta\}$ is countable and $\chi_\eta$
does not depend on $\omega_1\setminus J_\eta$. Therefore, the sets
$\{\chi_\zeta: \zeta\leq\eta\}$ and  $\{J_\zeta: \zeta\leq\eta\}$
satisfy (i)--(iv) at the step $\eta$.

For every $A\in \xi$, we put $F_A=\bar{\{\chi_\eta \colon \eta \in
A\}}^{\,\T^{G_\xi}}$ and $\mathcal{C}=\{F_A\colon A\in \xi\}$. It
follows from $\chi_\eta\in K$ for all $\eta<\omega_1$ and the
compactness of $K$ that $F_A\subset K$, for each $A\in\xi$. Since
$\mathcal{C}$ is a family of closed subsets of the compact space
$\T^{G_\xi}$ with the finite intersection property, $\mathcal{C}$
has non-empty intersection. Let $\rho$ be a point in $\bigcap
\{F_A\colon A\in \xi\}$. Clearly, $\rho\in K$, so $\rho$ is continuous.
Let $J_\rho$ be a countable subset of $\omega_1$ such that
$\rho$ does not depend on $\omega_1\setminus J_\rho$.

Since $\Sigma$ is a dense subgroup of both $G_\xi$ and $\Pi=
P\Z_2^{\omega_1}$, the characters $\rho$ and $\chi_\eta$ admit
continuous extensions $\bar{\rho}\colon\Pi\to \T$ and
$\bar{\chi}_\eta \colon\Pi\to \T$, for each $\eta<\omega_1$. Again,
the density of $G_\xi$ in $\Pi$ implies that $\bar\rho$ does not
depend on $\omega_1\setminus J_\rho$ and $\bar{\chi}_\eta$ does
not depend on $\omega_1\setminus J_\eta$.

Denote by $\mathbf 1$ the element of $\Pi$ all of whose coordinates
are equal to $1$. For every $\eta<\omega_1$, let $H_\eta=\{x\in
\Pi\colon x(\eta)=0\}$ and take a character $\psi_{\eta}$ on
$\Pi=\langle H_\eta, \mathbf{1}\rangle$ defined by
$\psi_\eta(x) = \bar{\chi}_\eta(x)$ and $\psi_\eta(x+\mathbf{1})=
\bar{\chi}_\eta(x)+\bar{\rho}(\mathbf{1})$, for all $x\in H_\eta$.
Since $\eta\in J_\eta$, we have $U(J_\eta) \subset H_\eta$. It then
follows from the above definition that $\psi_\eta$ does not depend
on $\omega_1\setminus J_\eta$.\smallskip

\textbf{Claim 1.} Put $T_\eta=\bigcup_{\zeta < \eta}J_\zeta$. For
every $\eta<\omega_1$, the character $\psi_\eta$ depends on
$\omega_1\setminus T_\eta$.\smallskip

\textit{Proof of Claim~1.} Indeed, by (iv) of the recursive
construction, $\chi_\eta$ depends on $\omega_1\setminus (\{\eta\}\cup
T_\eta)$. Hence there exists $x\in G_\xi\cap U(\{\eta\}\cup T_\eta)$
such that $\chi_\eta(x)\neq 1$. Then $x\in H_\eta$ and
$\psi_\eta(x)= \bar{\chi}_\eta(x)=\chi_\eta(x)\neq 1$, and we see
that $\psi_\eta$ depends on $\omega_1\setminus T_\eta$.\smallskip

\textbf{Claim 2.} For all $\alpha <\omega_1$, $\bar{\rho}\in
\bar{\{\psi_\eta \colon \alpha\leq \eta<\omega_1\}}^{\,\T^\Pi}$.
\smallskip

\textit{Proof of Claim~2.} Fix $\alpha <\omega_1$ and take
$\{g_1,\ldots,g_n\}\subset \Pi$. We can assume that there exists
$m\leq n$ such that $\{g_1,\ldots,g_m\}\subset G_\xi$ and
$\{g_{m+1},\ldots,g_n\}\subset \Pi\setminus G_\xi$. Then
$\{g_1,\ldots,g_m,g_{m+1}+\mathbf{1},\ldots,g_n+\mathbf{1}\}\subset
G_\xi$. Let $A\hskip-1pt=\omega_1\setminus\bigcup_{i\leq
m}\supp(g_i)$, $B\hskip-1pt=\bigcap_{m<k\leq n} \supp(g_k)$, and
$C=A\cap B\cap [\alpha,\omega_1)$. It follows from our choice of
$g_1,\ldots, g_n$ that $C\in \xi$. So, $\rho\in F_C$. Take $\eta\in
C$ such that $\rho (g_i)=\chi_\eta(g_i)$ whenever $1\leq i\leq m$
and $\rho (g_k+\mathbf{1})=\chi_\eta(g_k+\mathbf{1})$ whenever
$m<k\leq n$. If $1\leq i\leq m$, then $g_i(\eta)=0$ because $\eta\in
A$. So, $g_i\in H_\eta$ and $\psi_\eta(g_i)=\bar{\chi}_\eta(g_i)=
\chi_\eta(g_i)=\rho(g_i)= \bar{\rho}(g_i)$. If $m<k\leq n$, then
$(g_k+\mathbf{1})(\eta)=0$ because $\eta\in B$. So,
$(g_k+\mathbf{1})\in H_\eta$ and $\psi_\eta(g_k+\mathbf{1})=
\bar{\chi}_\eta(g_k+\mathbf{1})=\chi_\eta(g_k+\mathbf{1})=
\rho(g_k+\mathbf{1})= \bar{\rho}(g_k+\mathbf{1})$. Since $\psi_\eta$
and $\bar{\rho}$ are homomorphisms and $\psi_\eta(\mathbf{1})=
\bar{\rho}(\mathbf{1})$, we see that $\psi_\eta(g_k)=
\bar{\rho}(g_k)$. Therefore, $\bar{\rho}\in \bar{\{\psi_\eta \colon
\alpha\leq \eta<\omega_1\}} ^{\,\T^\Pi}$. This completes the proof
of Claim~2.\smallskip

Now, we have a character $\bar\rho\in\Pi^\wedge$, a family of
characters $\{\psi_\eta\colon \eta<\omega_1\}\subset \Pi^\wedge$,
and a family $\{J_\eta: \eta<\omega_1\}$ of countable subsets of
$\omega_1$. If follows from our definition of the characters
$\psi_\eta$'s and the above conditions (i), (ii), and (v) that
$\{\psi_\eta\colon \eta<\omega_1\}$ and $\{J_\eta: \eta<\omega_1\}$
satisfy (a)--(c) of Lemma~\ref{Le:MT} (with $\psi_\eta$'s in place
of $\chi_\eta$'s). Since $\bar{\rho}\in \bigcap_{\alpha<\omega_1}
\bar{\{\psi_\eta \colon\alpha\leq \eta<\omega_1\}}^{\,\T^\Pi}$, we are
in position to use Proposition~\ref{construction} to obtain a contradiction
with the fact that $\bar{\rho}$ is continuous. This contradiction shows that
the compact set $K$ depends at most on countably many coordinates.
The reflexivity of $G_\xi$ now follows from Theorem~\ref{duality}.
\end{proof}

To finish this section, we extend the conclusion of Proposition~\ref{main1}
to subgroups $G_\xi$ of the group $P\Z_2^\tau$, for any uncountable
cardinal $\tau$.

\begin{theorem}\label{FinSec}
Let $\tau>\omega$ be a cardinal and $\xi$ an ultrafilter on $\tau$
containing the complements to countable sets. Then subgroup $G_\xi$
of $P\Z_2^\tau$ is reflexive.
\end{theorem}

\begin{proof}
According to Theorem~\ref{Reduc} it suffices to verify that the subgroup
$\pi_J(G_\xi)$ of the group $P\Z_2^J$ is reflexive, for every set $J\subset
\tau$ satisfying $|J|\leq\aleph_1$. Let us consider two possible cases.

Case~1. $J\notin\xi$. Then the definition of $G_\xi$ implies that
$\pi_J(G_\xi)=\Z_2^J$, so the reflexivity of $\pi_J(G_\xi)$ is immediate
from Theorem~\ref{cor:main}.

Case~2. $J\in\xi$. Put $\eta=\{J\cap A: A\in\xi\}$. Then $\eta$ is an ultrafilter
on $J$ containing the complements to countable sets. Further, the definition
of $G_\xi$ implies that the projection $\pi_J(G_\xi)$ of $G_\xi$ coincides with
the subgroup $G_\eta$ of $P\Z_2^J$. Identifying $J$ and $\omega_1$ and
applying Proposition~\ref{main1}, we see that the group $\pi_J(G_\xi)$ is again
reflexive.
\end{proof}

%%%%%%%%%%%%%%%%%%%%%%%%%%%%%%%%%%
\section{Non-reflexive $P$-groups}\label{NonRef}%%%%%%%%%
We would like to trace the border between reflexivity and
non-reflexivity for $P$-groups $G$ such that $\Sigma{D}\subset G
\subset PD$, where $D=\prod_{i\in I}D_i$ is a product of discrete
groups. Recall that by an old result of Leptin in \cite{lept} (see also
\cite[Example~17.11]{bana91}), the subgroup
\[
\sigma\Z_2^{\omega_1}=\{x\in\Z_2^{\omega_1} \colon \supp(x)
\mbox{ is finite\,}\}
\]
of $P\Z_2^{\omega_1}$ is not reflexive. We now extend this fact to
some dense subgroups of the groups of the form $PD$.

Let $G$ be a subgroup of $PD$ containing $\Sigma{D}$. For each
$i\in I$, let $\pi_i\colon G \to D_i$ be the projection, $\pi_i(x)=x(i)$.
For a set $J\subset I$, we also put $F_J=
\overline{\{ \varphi \circ \pi_i \colon i \in J,\  \varphi \in
(D_i)^\wedge \}}^{\,\T^G} \bigcup \{\mathbf{1}\}$, where $\mathbf{1}$
is the identity of $G^\wedge$.

\begin{lemma}\label{countclos1}
For every $J\in [I]^{\leq\omega}$, $F_J\subset G^\wedge$.
\end{lemma}

\begin{proof}
Suppose that $J\in [I]^{\leq\omega}$ and take any $\rho\in F_J$.
Then $U(J)$ is an open set in $PD$ containing $0_G$ such that
$\rho(U(J)\cap G)=\{1\}$. Therefore $\rho$ is continuous.
\end{proof}

\begin{lemma}\label{cont}
Let $J$ be a nonempty subset of $I$ and $\psi\in F_J$,
$\psi\neq\mathbf{1}$. Then $\psi$ is continuous as a character
on $G$ if and only if there exists $x\in \Sigma{D}$ such that
$\psi(x)\neq 1$.
\end{lemma}

\begin{proof}
Suppose that $\psi$ is continuous. Since $\psi\neq\mathbf{1}$
and $\Sigma{D}$ is dense in $G$, there is a point $x\in \Sigma{D}$
such that $\psi(x)\neq 1$.

Conversely, take $x\in \Sigma{D}$ such that $\psi(x)\neq 1$ and
write $F_J=F_{J\cap \supp(x)}\bigcup F_{J\setminus \supp (x)}$.
Since $(\chi \circ \pi_i)(x)=1$ for all $i\in J\setminus \supp(x)$
and all $\chi \in (D_i)^\wedge$, we have that $\psi \notin
F_{J\setminus \supp (x)}$. Then $\psi \in F_{J\cap \supp(x)}$. It
now follows from Lemma~\ref{countclos1} that $\psi$ is continuous.
\end{proof}

\begin{proposition}\label{sigmanot}
Let $D=\prod_{i\in I}D_i$ be a product of nontrivial discrete
Abelian groups. Then the subgroup $\Sigma=\Sigma{D}$ of
$PD$ is not reflexive provided that $|I|>\omega$. Furthermore,
the bidual group $\Sigma^{\wedge \,\wedge}$ is discrete.
\end{proposition}

\begin{proof}
Given $\psi\in F_I$, $\psi\neq\mathbf{1}$, there exists $x\in
\Sigma$ such that $\psi(x)\neq 1$. It now follows from
Lemma~\ref{cont} that $F_I\subset\Sigma^\wedge$. The set $F_I$ is
compact as a closed subset of $\T^{\Sigma}$. Since $F_I$ does not
depend  on countably many coordinates (it actually  depends on every
index $i\in I$), the group $\Sigma$ is not reflexive by
Theorem~\ref{duality}. It is easy to see that $F_I$ generates a dense
subgroup of $\Sigma^\wedge$, so $\left(F_I\right)
^\vartriangleright$ contains only the neutral element of
$G^{\wedge\wedge}$. Hence the bidual group $\Sigma^{\wedge
\,\wedge}$ is discrete.
\end{proof}

Sets of the form $F_J$, with $J\neq I$,  can also be used to show
that some subgroups larger than $\Sigma=\Sigma{D}$ are not reflexive.
This is done in Lemma \ref{NonRef1} for groups of the form
$G_L=\langle\Sigma, L\rangle$ with  $L\subset PD$ satisfying
$|I\setminus\bigcup _{x\in L}\supp(x)|\geq \omega_1$.

\begin{lemma}\label{comp}
Let $L$\ be a subset of $PD$ and $J=I\setminus\bigcup_{x\in L}
\supp(x)$. If $|J|>\omega$, then $F_J=\overline{\{\varphi \circ \pi_i
\colon i\in J,\ \varphi \in (D_i)^\wedge\}} ^{\,\T^{G_L}}\bigcup \{\mathbf{1}\}$
is a compact subset of $(G_L)^{\wedge}$.
\end{lemma}

\begin{proof}
Take $\psi \in F_J$, $\psi \neq\mathbf{1}$. Assume that $\psi (z)=1$
for all $z\in\Sigma$. Take $g\in G$ such that $\psi(g)\neq 1$ and write
$g=z+n_1 x_1+\ldots+n_k x_k$, where $z\in\Sigma$, $x_1,\ldots,x_k
\in L$, and $n_1,\ldots,n_k\in\Z$. Since $\psi (g)\neq 1$ and
$U=\{y\in\T^G: y(g)\neq 1\}$ is an open neighbourhood of $\psi$, there
are $j\in J$ and $\varphi\in(D_j)^\wedge$ such that $(\varphi \circ\pi_j)
(g)\in U$. Then $(\varphi\circ \pi_j)(g)\neq 1$ which is a contradiction
because
\[
(\varphi \circ \pi_j)(g)=(\varphi \circ \pi_j)(n_1x_1+\cdots+n_k x_k)=
\varphi (n_1 x_1(j)+\cdots+ n_k x_k(j))=\varphi (0_G)=1.
\]
So, there is $z\in \Sigma$ such that $\psi (z)\neq 1$. Now the
continuity of $\psi$ follows from Lemma~\ref{cont}.
\end{proof}

\begin{theorem}\label{NonRef1}
Let $L$ be a subset of $PD$ such that the set $J=I\setminus
\bigcup_{x\in L}\supp(x)$ is uncountable. Then the subgroup
$G=\langle \Sigma, L\rangle$ of $PD$ is not reflexive.
\end{theorem}

\begin{proof}
Put $F_J=\overline{\{ \varphi \circ \pi_i \colon i\in J,\ \varphi \in
(D_i)^\wedge\}}^{\,\T^G}$. It follows from Lemma~\ref{comp}
that $F_J\subset G^\wedge$. Since $F_J$ depends on every
coordinate $i\in J$, we conclude that $G$ is not reflexive in view
of Theorem~\ref{duality}.
\end{proof}

The preceding results make use of  sets of the form $F_J$ to see
that some subgroups $G$ with $\Sigma{D}\subseteq G\subset PD$
are not reflexive. Sets of this sort were already used by Leptin \cite{lept}
to evidence the nonreflexivity of $\sigma\Z_2^{\omega_1}$. We see
next that  $F_J$ may not be contained in $G^\wedge$, for some
$G=\langle \Sigma{D}, a\rangle $ and $J\subset I$. This makes necessary
a  different approach to show, as we do in Theorem \ref{GLnot}, that
these groups are not reflexive.

Let $\Pi=PD$, $\Sigma=\Sigma{D},$ and $a$ be a point in
$\Pi \backslash \Sigma$ such that $na\in\Sigma$ for some integer
$n>1$. We put $G=\langle \Sigma,a\rangle$, $J=\supp (a)$, and
$F_J=\overline{ \{ \varphi \circ \pi_i\colon i\in J,\ \varphi \in
(D_i)^\wedge\}}^{\,\T^{G}}$. Note that
$G=\{g+ka\colon g\in \Sigma,\ 0\leq k <n\}$.

\begin{lemma}\label{Le:FIdontwork}
If $G=\langle \Sigma,a\rangle$, then $F_J$ is not contained in
$G^\wedge$.
\end{lemma}

\begin{proof}
Since $na\in\Sigma$, there exists a countable set $C\subset I$ such
that $na(i)=0_i$ for each $i\in I\setminus C$. Therefore, we can
find a divisor $m$ of $n$ with $m>1$ and an uncountable set
$A\subset J\setminus C$ such that the order of $a(i)$ equals $m$ for
each $i\in A$.

Let us define $\psi$ in Hom($G,\T$) as follows: $\psi(g+ la) =
e^{(2\pi l/m) i}$ for all $g\in \Sigma$ and $l\in\Z$. Since
$\psi(g)=1$ for each $g\in\Sigma$ and $\psi(l_1 a)=\psi(l_2 a)$
whenever $m$ divides $l_1-l_2$, our definition of $\psi$ is correct.
We now show that $\psi \in F_A$. For this, let $g_k$ be points in
$\Sigma$, $l_k$ be integers, and $V_k$ open sets in $\T$ such that
$\psi(g_k+l_k a)\in V_k$ for all $k\leq N$, where $N\in\N$. Since
$|\bigcup_{0\leq k\leq N} \supp(g_k)|\leq \omega$ and $|A|\geq
\omega_1$, we can choose $j \in A\backslash \bigcup_{0\leq k\leq N}
\supp(g_k)$. Since the order of $a(j)$ equals $m$, there is a
character $\rho\in (D_j)^\wedge$ such that $\rho(a(j))=e^{(2\pi/m)i}$.
Now it is easy to see that $(\rho \circ \pi_j)(g_k + l_ka)=
\psi(g_k+l_k a) \in V_k$ whenever $0\leq k\leq N$. Therefore
$\psi \in F_A\subset F_J$. The discontinuity of $\psi$ follows from
Lemma~\ref{cont} because $\psi\neq\mathbf{1}$ and $\psi(g)=1$ for
all $g\in \Sigma$. This proves that $\psi\in F_J\setminus
G^\wedge\neq\emptyset$.
\end{proof}

By Lemma \ref{Le:FIdontwork}, the argument used in
Proposition~\ref{sigmanot} for $\Sigma{D}$ does not work for
$\langle\Sigma,\overline{1}\rangle$, where $\Sigma=
\Sigma\Z_2^{\omega_1}$ and $\overline{1}$ is the element of
$\Z_2^{\omega_1}$ all of whose coordinates are equal to $1$.
However, we show in Theorem~\ref{GLnot} that the group
$\langle\Sigma,\overline{1}\rangle$ is not reflexive either.

\begin{lemma}\label{NINS}
Let $\Pi=P\Z_2^{\omega_1}$ and $\Sigma=\Sigma\Pi$. Suppose
that a subgroup $L$ of $\Pi$ has the property that for each
$\alpha<\omega_1$, there exists an ordinal $\beta(\alpha)$ such
that $\alpha<\beta(\alpha)<\omega_1$ and the restriction to $L$ of
the projection $p_{J(\alpha)}\colon \Pi \to\Z_2^{J(\alpha)}$ is injective,
where $J(\alpha)=\beta(\alpha)\setminus\alpha$. Then the subgroup
$G=\Sigma+L$ of $\Pi$ is not reflexive.
\end{lemma}

\begin{proof}
For every $\alpha<\omega_1$, let $\gamma(\alpha)=\beta(\alpha)+1$
and decompose $\Z_2^{\gamma(\alpha)}=H_\alpha \oplus L_\alpha$,
with $L_\alpha=p_{\gamma(\alpha)}(L)$, where $p_{\gamma(\alpha)}
\colon \Pi\to \Z_2^{\gamma(\alpha)}$ is the projection. Observe that
we can (and so we will) choose $H_\alpha$ to satisfy the following
condition:
\begin{equation}\label{1}
\left\{ x\in \Z_2^{\gamma(\alpha)}\colon
x(\eta)=0 \mbox{ for all } \eta\in J(\alpha)\right\}\subset H_\alpha.
\end{equation}
We also choose, for every $\alpha<\omega_1$, a homomorphism
$\psi_\alpha\colon \Z_2^{\gamma(\alpha)}\to \Z_2$ such that
$\psi_\alpha(L_\alpha)=\{0\}$ and $\psi_\alpha\hskip-4pt
\upharpoonright\hskip-4pt_{H_\alpha}=\pi_{\beta(\alpha)}
\hskip-4pt\upharpoonright\hskip-4pt_{H_\alpha}$, where
$\pi_{\beta(\alpha)}$ is the projection of $\Z_2^{\gamma(\alpha)}$
to the $\beta(\alpha)$th factor $(\Z_2)_{\beta(\alpha)}$.
Let us put $\chi_\alpha=\psi_{\alpha}\circ p_{\gamma(\alpha)}$.
Since both homomorphisms $\psi_\alpha$ and $p_{\gamma(\alpha)}$
are continuous, so is $\chi_\alpha$.

For every $\alpha\le\omega_1$, let $E_\alpha=\overline{
\{\chi_\sigma\colon\sigma<\alpha\}}^{\,\Z_2^G}$. We claim that
$E_{\omega_1}\subset G^\wedge$. Indeed, let $\rho \in
E_{\omega_1}$. If $\rho \neq \mathbf{0}$, there is $g\in G$ with
$\rho(g)=1$. Take $g_1\in \Sigma$ and $g_2\in L$ such that
$g=g_1+g_2$. If $\rho(g_2)=1$, then there exists $\alpha<\omega_1$
such that $\chi_\alpha(g_2)=1$ which is not possible because
$\chi_\alpha(g_2)=\psi_{\alpha}(p_{\gamma(\alpha)}(g_2))=0$. We thus
have $\rho(g_1)=1$. Now, mimicking the proof of Lemma \ref{cont}, we
may take $\alpha<\omega_1$ such that $\supp(g_1)\subset \alpha$, and
write
\[
E_{\omega_1}=E_\alpha \bigcup \overline{\{ \chi_\sigma \colon
\sigma\geq\alpha \}}^{\,\Z_2^G}.
\]
Take any $\sigma\geq \alpha$. It follows from (1) that
$p_{\gamma(\sigma)}(g_1)\in H_\sigma$ and hence
$\chi_\sigma(g_1)\hskip-1pt=\hskip-1pt\psi_{\sigma}
(p_{\gamma(\sigma)}(g_1))\hskip-1pt=\pi_{\beta(\sigma)}
(p_{\gamma(\sigma)}(g_1))=0$. Since $\sigma\geq\alpha$ is arbitrary,
we conclude that $\rho\notin\overline {\{\chi_\sigma\colon
\sigma\geq\alpha \}}^{\,\Z_2^G}$ and, therefore, $\rho\in E_\alpha$.
We have thus that $E_{\omega_1}=\bigcup_{\alpha<\omega_1} E_\alpha$.
Since all elements of $E_\alpha$ are continuous (they all do not
depend on $\omega_1\setminus \gamma(\alpha)$), we see that
$E_{\omega_1}$ is a compact subset of $G^\wedge$.

Given $\alpha<\omega_1$, we define $b_\alpha\in\Sigma$
by $b_\alpha (\gamma)= 0$ if $\gamma \neq \beta(\alpha)$
and $b_\alpha(\beta(\alpha))=1$. This implies that
$p_{\gamma(\alpha)}(b_\alpha)\in H_\alpha$ and, therefore,
$\chi_\alpha(b_\alpha)=\psi_\alpha(p_{\gamma(\alpha)}(b_\alpha))=
\pi_{\beta(\alpha)}(p_{\gamma(\alpha)}(b_\alpha))=
b_\alpha(\beta(\alpha))=1$. It follows that $\chi_\alpha$ depends on
the index $\beta(\alpha)$ which in its turn implies that $E_{\omega_1}$
depends on uncountably many coordinates. By Theorem~\ref{duality},
the group $G$ is not reflexive.
\end{proof}

\begin{theorem}\label{GLnot}
Let $\tau$ be an uncountable cardinal, $\Pi=P\Z_2^\tau$, and
$\Sigma=\Sigma\Pi$. Then, for every countable subgroup $L$ of $\Pi$,
the group $G=\Sigma+L\subset\Pi$ is not reflexive.
\end{theorem}

\begin{proof}
Let $L_0=L\cap\Sigma$. Since every subgroup of the boolean group
$L$ is a direct summand, there exists a subgroup $L_1$ of $L$ such
that $L=L_0\oplus L_1$. Since $L_0\subset\Sigma$, we see that
$G=\Sigma+L=\Sigma+L_0+L_1=\Sigma+L_1$. It follows from our
definition of $L_1$ that the intersection $\Sigma\cap L_1$ is trivial,
so the set $\supp(x)$ is uncountable, for each $x\in L_1$ distinct
from $0_G$.

Take a subset $J$ of $\tau$ such that $|J|=\aleph_1$ and
$|J\cap\supp(x)|=\aleph_1$ for each $x\in L_1$, $x\neq 0_G$. By
Theorem~\ref{Reduc}, it suffices to show that the subgroup
$\pi_J(G)$ of $P\Z_2^J$ is not reflexive. Since $\pi_J(G)=
\pi_J(\Sigma)+\pi_J(L_1)$ and $\pi_J(\Sigma)=\Sigma\Z_2^J$,
we can assume without loss of generality that $\tau=J=\omega_1$.
Hence $G$ is a subgroup of $\Pi=P\Z_2^{\omega_1}$.

In view of Lemma~\ref{NINS} it suffices to show that for every
$\alpha<\omega_1$, there is a countable ordinal
$\beta(\alpha)>\alpha$ such that the restriction to $L_1$ of the
projection $p_{[\alpha,\beta(\alpha))}$ is injective. Given an
ordinal $\alpha<\omega_1$, we take an element $\beta(x)\in
\supp(x)\setminus\alpha$, for each $x\in L_1\setminus\{0_G\}$,
and put
\[
B=\{\beta(x): x\in L_1,\  x\neq 0_G\}.
\]
Then $|B|\leq |L_1|\leq\omega$, so there exists a countable ordinal
$\beta(\alpha)>\alpha$ such that $B\subseteq \beta(\alpha)$. Put
$J(\alpha)=\beta(\alpha)\setminus\alpha$. It is clear that the
restriction to $L_1$ of the projection $p_{J(\alpha)}$ is one-to-one
since $B\subset J(\alpha)$ and $B$ intersects the set $\supp(x)
\setminus\alpha$, for each $x\in L_1$ distinct from $0_G$. This
finishes the proof of the theorem.
\end{proof}

%%%%%%%%%%%%%%%%%%%%%%%%%%%%%%%%%%%%%
\section{Problem section}\label{Sec:Problems}
Here we present several problems whose solutions can substantially improve
our understanding of the duality theory for $P$-groups. The first of them
arises in an attempt to extend Theorem~\ref{cor:main} to arbitrary products
of reflexive $P$-groups:

\begin{problem}\label{QM0}
Let $\Pi=\prod_{i\in I}G_i$ be the product of a family of reflexive $P$-groups.
Is the group $P\Pi$ reflexive?
\end{problem}

According to Proposition~\ref{Pro:Ge} it suffices to consider the case when
the index set $I$ in the above problem is countable. One can try to prove
(or refute) a more general form of the above problem inspired by
Proposition~\ref{Cor:Comp}:

\begin{problem}\label{New}
Let $G$ be a reflexive topological group. Is the group $PG$ then reflexive?
\end{problem}

A direct verification shows that every reflexive $P$-group constructed so far
contains a discrete (hence closed) subgroup of cardinality $2^\omega$. This
explains the origin of the following problem:

\begin{problem}\label{LPG}
Does there exist a nondiscrete reflexive Lindel\"of $P$-group?
\end{problem}

In the next problem we pretend to generalize Theorem~\ref{FinSec}:

\begin{problem}\label{QM1}
Let $\tau$ be an uncountable cardinal and $G$ be a subgroup of
$\Pi=P\Z_2^\tau$ such that $|\Pi/G|<\omega$ (or
$|\Pi/G|\leq\omega$). Is $G$ reflexive? What if, additionally,
$G$ contains $\Sigma\Pi$?
\end{problem}

We do not know whether Theorem~\ref{GLnot} extends to bigger
subgroups of $P\Z_2^{\omega_1}$:

\begin{problem}\label{QM2}
Is is true that the subgroup $\Sigma+L$ of $\Pi=P\Z_2^{\omega_1}$
fails to be reflexive, for any subgroup $L$ of $\Pi$ satisfying
$|L|\leq\aleph_1$?
\end{problem}

\def\cprime{$'$} \def\cprime{$'$}
  \def\polhk#1{\setbox0=\hbox{#1}{\ooalign{\hidewidth
  \lower1.5ex\hbox{`}\hidewidth\crcr\unhbox0}}}
  \def\polhk#1{\setbox0=\hbox{#1}{\ooalign{\hidewidth
  \lower1.5ex\hbox{`}\hidewidth\crcr\unhbox0}}} \def\cprime{$'$}
  \def\cprime{$'$} \def\cprime{$'$}
  \def\polhk#1{\setbox0=\hbox{#1}{\ooalign{\hidewidth
  \lower1.5ex\hbox{`}\hidewidth\crcr\unhbox0}}}
  \def\polhk#1{\setbox0=\hbox{#1}{\ooalign{\hidewidth
  \lower1.5ex\hbox{`}\hidewidth\crcr\unhbox0}}}
\providecommand{\bysame}{\leavevmode\hbox
to3em{\hrulefill}\thinspace}
\providecommand{\MR}{\relax\ifhmode\unskip\space\fi MR }
\providecommand{\MRhref}[2]{%
  \href{http://www.ams.org/mathscinet-getitem?mr=#1}{#2}
} \providecommand{\href}[2]{#2}

%\bibliography{C:/Users/Admin/Documents/BIBLIOGRAFIAS/bibrep}
%C:/Users/Admin/Documents/BIBLIOGRAFIAS/bibbohr,C:/Users/Admin/Documents/
%BIBLIOGRAFIAS/bibsidon}
%\bibliography{E:/DOCS/BIBLIOGRAFIAS/bibrep}
%\bibliography{F:/bibliografias/bibrep}

\end{document}